\documentclass[11pt]{scrartcl}
\usepackage[final]{allergy}
\usepackage{unixode}
\usepackage{helvetic}
\ifxetex
\RequirePackage{fontspec}
\else
\usepackage{inconsolata}
\fi

\usepackage{booktabs}
\usepackage{enumitem}

\usepackage{float}
\usepackage[hang]{subfigure}

\usepackage{authblk}

\usepackage{pythonhighlight}

\newcommand{\Sfrac}[2]{{\textstyle{\frac{#1}{#2}}}}

\newcommand{\R}{\mathbb{R}}
\newcommand{\NN}{\mathbb{N}}

\newcommand{\Lop}{\mathcal{L}}
\newcommand{\Aop}{\mathcal{A}}
\newcommand{\Hop}{\mathcal{H}}
\newcommand{\DD}{\mathcal{D}}
\newcommand{\SN}{\mathcal{S}_N}

\newcommand*\dd{\mathrm{d}}
\newcommand*\ii{\mathrm{i}}

\DeclareMathOperator{\sech}{sech}

\newcommand*\email[1]{\htmladdnormallink{\texttt{#1}}{mailto:#1}}

\author[a]{Henrik Kalisch\thanks{\email{henrik.kalisch@math.uib.no}}}
\author[a]{Daulet Moldabayev\thanks{\email{daulet.moldabayev@math.uib.no}}}
\author[b,c]{Olivier Verdier\thanks{\email{olivier.verdier@hib.no}} }
\affil[a]{Department of Mathematics, University of Bergen }
\affil[b]{Department of Mathematics and Statistics, University of Umeå, Sweden } 

\affil[c]{Department of Computing Mathematics and Physics, Høgskolen i Bergen }

\date{}

\title{Numerical Study of Nonlinear Dispersive Wave Models with SpecTraVVave}

\graphicspath{{./figures/}}

\begin{document}

\maketitle


\begin{abstract}
In nonlinear dispersive evolution equations, the competing effects of nonlinearity and dispersion 
make a number of interesting phenomena possible. In the current work, the focus is on the numerical
approximation of traveling-wave solutions of such equations.
We describe our efforts to write a dedicated \textsf{Python}
code which is able to compute traveling-wave solutions of nonlinear
dispersive equations of the general form
\begin{equation*}
u_t + [f(u)]_{x} + \Lop u_x = 0,
\end{equation*}
where $\Lop$ is a self-adjoint operator, and $f$ is a real-valued function with  $f(0)  = 0$.

The \textsf{SpectraVVave} code uses a continuation method coupled with a spectral projection
to compute approximations of steady symmetric solutions of this equation.
The code is used in a number of situations to gain an
understanding of traveling-wave solutions. The first case is the Whitham
equation, where numerical evidence points to the conclusion that the main 
bifurcation branch features three distinct points of interest, namely a turning point, 
a point of stability inversion, and a terminal point which corresponds to a cusped wave.

The second case is the so-called modified Benjamin--Ono equation where
the interaction of two solitary waves is investigated. It is found
that is possible for two solitary waves to interact in such a way
that the smaller wave is annihilated. The third case concerns
the Benjamin equation which features two competing dispersive operators.
In this case, it is found that bifurcation curves of periodic traveling-wave
solutions may cross and connect high up on the branch in the nonlinear regime.

\end{abstract}

\section{Introduction.}
%
%
%
%

This paper is concerned with traveling wave solutions for a class of nonlinear dispersive
equations of the form
\begin{equation}\label{nonlocal}
u_t + [f(u)]_x + \Lop u_x = 0,
\end{equation}
where $\Lop$ is a self-adjoint operator, and $f$ is a real-valued function 
with  $f(0)=0$ and which satisfies certain growth conditions.
Equations of this form arise routinely in the study of wave problems in fluid mechanics
and mony other contexts.
A prototype of such an equation is the KdV equation that appears
if $\Lop= I + \frac{1}{6} \partial_x^2$ and $f(u) = \frac{3}{4} u^2$.
In the current work, the operator $\Lop$ is considered to be given as a Fourier multiplier operator,
such as for instance in the Benjamin--Ono equation, which arises in the study of
interfacial waves. In this case, the Fourier multiplier operator is given by
$\Lop =  I - \Hop \partial_x$, where the Hilbert transform
$\Hop$ is defined as
\begin{equation}
\Hop u(x) = \frac{1}{\pi} \, \mathrm{p.v.} \int_{-\infty}^{\infty} \frac{u(x-y)}{y} \, \dd y, 
\qquad \widehat{\Hop u}(k) = -\ii \operatorname{sgn}(k)\widehat{u}(k).
\end{equation}
We also study in detail traveling wave solutions of the Whitham equation, which appears when 
$\Lop$ is given by convolution with the integral kernel $K_{h_0}$ in the form
\begin{equation}\label{symbol}
\Lop u(x) = \int_{-\infty}^{\infty} K_{h_0} (y) u(x-y) \, \dd y,
\qquad \widehat{K_{h_0}} (k) 
= {\textstyle{ \sqrt{\frac{g\tanh (h_0 k)}{k}}}}, 
\end{equation}
and $f$ is the same function as in the KdV equation.

The particular form of equation~\eqref{nonlocal} exhibits the competing effects of dispersion and nonlinearity, 
which gives rise to a host of interesting phenomena. The most well known special phenomenon is the existence
of solitary waves and of periodic traveling waves containing higher Fourier modes.
Indeed, note that in the purely dispersive model $u_t + \Lop u_x = 0$, 
the only possible permanent progressive waves are simple sinusoidal waves, while
in the nonlinear model~\eqref{nonlocal} higher Fourier modes must be considered
to obtain solutions.

The order of the operator $\mathcal{L}$ appearing in~\eqref{nonlocal} 
has a major effect on the types of solutions that may be found. 
A higher-order operator, such as in the Korteweg--de Vries equation, 
acts as a smoothing operator because of its effect of spreading different frequency components
out due to a strongly varying phase speed~\cite{Kato}. Lower-order operators 
such as the operator $K_{h_0}$ in~\eqref{symbol} appearing in the Whitham equation
may allow solutions to develop singularities, such as derivative blow-up (see~\cite{Hur2015, HurTao14})
and formation of cusps (see~\cite{EW2016}).

On the other hand, highly nonlinear functions $f(u)$ may lead to $L^{\infty}$-blow-up.
For instance, the generalized KdV equation which is written in normalized form as
\begin{equation}
\label{kdv}
u_t + u^{p}u_x + u_x + u_{xxx} = 0,
\end{equation}
features global existence of solutions for $p=1,2,3$, 
but the solutions blow-up in the critical case $p=4$ (the case $p \ge 5$ is open).
In the case of the generalized Benjamin--Ono equation
\begin{equation*}
u_t + u^{p}u_x + u_x - \mathcal{H} u_{xx} = 0,
\end{equation*}
where $\mathcal{H}$ is the Hilbert transform, numerical evidence points to singularity formation 
for $p > 2$~\cite{BoKa}, but no proofs are available at this time.

In order to study different phenomena related to equations of the form~\eqref{nonlocal} 
and their traveling wave solutions, a \emph{Python}-based solver package \textsf{SpectraVVave} 
was developed by the authors~\cite{github}. 
The general idea behind the solver is to use a numerical continuation method~\cite{Keller} 
implemented with a pseudo-spectral algorithm. 
Similar previous projects include \textsf{AUTO}~\cite{Doed97} and \textsf{Wavetrain}~\cite{Sherratt}.
\textsf{AUTO} is written in \emph{C}, whereas \textsf{Wavetrain} is written in \emph{Fortran}. 
Both programs are efficient and very general, as they are able to cover a wide range of problems
involving bifurcation analyses. However, from a user's perspective, 
such a generality coupled with low level programming languages
may lead to some difficulty for users of these programs to utilize them efficiently.

\textsf{SpectraVVave} is designed to provide researchers with a simple yet effective tool 
for investigating problems on traveling waves. The package is flexible, 
and its functionality can be easily expanded. 
The availability of the \textsf{IPython} notebook \cite{IPython} makes the solver very interactive,
so that it should be easier for new users to get started.

In order to maximize ease of use, \textsf{SpectraVVave} was designed
to find even solutions of~\eqref{nonlocal}.
Symmetry of steady solutions can be proved for some of the models
in the form~\eqref{nonlocal}, but nor for all~\cite{CB}.
Some of these equations also admit non-smooth solutions,
for instance as termination points of a bifurcation branch.
This happens for exmple for the Whitham equation, which features bifurcation
curves which terminate in a solution with a cusp~\cite{EW2016}.
One of the goals of the present paper is to investigate the
precise nature of the termination of the bifurcation curve.



The content of the paper is structured as follows. A mathematical description of the numerical method of 
\textsf{SpectraVVave} is given in \autoref{sec:spectral}. \autoref{sec:experiments} presents results of different experiments carried out with the package. Concluding remarks are given in \autoref{sec:conclusion}. A method for finding initial guesses for the solver is described in \autoref{sec:stokes}. \autoref{sec:implementation} contains a schematic of program and a description of its workflow.

\section{Spectral scheme and construction of nonlinear system.}
\label{sec:spectral}
\subsection{Cosine collocation method.}
To compute traveling wave solutions to the equation~\eqref{nonlocal} the following ansatz is employed:
\begin{equation}\label{ansatz}
u(x,t) = \phi(x-ct).
\end{equation}
Thus, the equation takes the form
\begin{equation}\label{phi-nonlocal}
\phi' + \left[f(\phi) \right]' + \Lop \phi' = 0,
\end{equation}
which can be integrated to give
\begin{equation}\label{int_gen1}
-c \phi + f(\phi)  + \Lop \phi = B. 
\end{equation}
The constant $B$ is a priori undetermined. One may set the $B$ equal to 
zero as a way of normalizing the solutions. Another option is to impose
an additional condition, for example that the integral of $\phi$ over
one wavelength be zero. In this case, $B$ will be found along with the solution $\phi$.

We consider $\Lop$ as a Fourier multiplier operator with symbol
$\alpha(k)$. We also assume that $f$ is at least twice differentiable,
and we have $f(0)=0$, $f'(0)=0$ and $f''(0)=0$.
When computing traveling wave solutions we focus on even periodic solutions. 
While it can be proved in ome cases that solutions of~\eqref{int_gen1}
must be even, this is not known for a general operator $\Lop$.
Nevertheless, we make this assumption here in order to
make the numerical procedure as uniform as possible.
For even periodic solutions, 
one may use a cosine collocation instead of a Fourier method. In particular, using
the cosine functions as basis elements automatically removes the inherent symmetries
due to reflective and translational symmetry. 
Moreover, the number of unknowns is reduced by a factor of $2$, 
and the problem of the asymmetric arrangement of nodes in the FFT
is circumvented. Of course, all these problems could also be dealt with a collocation
method based on the Fourier basis, but the cosine basis does all of the above automatically.
In addition,  the \emph{Python} cosine transform is based on an integrated algorithm, which
relies on an optimized version of the discrete cosine transform (DCT).
	
The following description of computation scheme was presented in detail in~\cite{EK2}, 
but we will briefly repeat it here for consistency of the manuscript. For the purpose of clarity, 
we will refer to full wavelength $L$ of a solution as fundamental wavelength, 
and a half of fundamental wavelength will be called wavelength. 
Such a definition is required because the method computes a half of a solution profile, 
the other half is automatically constructed due to symmetry.
 
Traveling wave solutions to the equation~\eqref{int_gen1}
are to be computed in the form of a linear combination of cosine functions of different
wave-numbers, i.e., in the space
\begin{equation}
\SN = \mathrm{span}_{\R}\setc{\cos(lx)}{0 \leq l \leq N-1 }
.
\end{equation}
This is a subspace of $L^2(0,2\pi)$, and the collocation points 
$x_n = \pi \frac{2n -1}{2N}$ for $n = 1,\dots, N$ are used to discretize the domain.  
If the required fundamental wavelength of solutions is to be $L \neq 2\pi$, 
one can use a scaling on the $x$-variable. 
Defining the new variable
\begin{equation}
x' = \frac{L}{2\pi}x, 
\end{equation}
yields collocation points $x_n'$ and wavenumbers $\kappa_l$ defined by
\begin{equation}
x_n' = \frac{L}{2} \frac{2n -1}{2N}, \quad\quad \kappa_l = \frac{2\pi}{L}l. \label{scalings}
\end{equation}
We are seeking a function $\phi_N \in \SN$ that satisfies the equations
\begin{equation}
-c\phi_N(x_n') + f(\phi_N)(x_n') + \Lop^N\phi_N(x_n') = 0, \label{int_gen-disc}
\end{equation}
at the collocation points $x_n'$.
The operator $\Lop^N$ is the discrete form of the operator $\Lop$, 
and $\phi_N$ is the discrete cosine representation of $\phi$
which is given by
\begin{align}
\phi_N(x') &= \sum_{l = 0}^{N-1} \omega(\kappa_l)\Phi_N(\kappa_l)\cos(\kappa_{l}x'),\\
\omega(\kappa_l) &= 
\begin{cases}
\sqrt{1/N}, \quad \kappa_l = 0, \\
\sqrt{2/N}, \quad \kappa_l > 0,
\end{cases}                  \\
\Phi_N(\kappa_l) &= \omega(\kappa_l) \sum_{n = 1}^{N} \phi_N(x'_n) \cos(\kappa_l x'_n),
\end{align}
where $\kappa_l = 0,\frac{2\pi}{L},\ldots, \frac{2\pi}{L}(N-1)$ are the scaled wavenumbers,
and $\Phi_N(\cdot)$ are the discrete cosine coefficients.
As the equation~\eqref{int_gen-disc} is enforced at the collocation points $x_n'$, 
one may evaluate the term $\Lop^N\phi_N$ using the matrix $\Lop^N(i,j)$ 
defined by
\begin{align}
\Lop^N\phi_N(x'_i) &= \sum_{j=1}^{N}\Lop^N(i,j) \phi_N(x'_j),\\
\Lop^N(i,j) &= \sum_{l=0}^{N-1} \omega^2(\kappa_l) \alpha(\kappa_l) \cos(\kappa_l x'_i) \cos(\kappa_l x'_j),
\end{align}
where $\alpha(\cdot)$ is the Fourier multiplier function of the operator $\Lop$.


\subsection{Construction of nonlinear system.}

%
The equation~\eqref{int_gen-disc} enforced at $N$ collocation points 
yields a nonlinear system of $N$ equations in $N$ unknowns, which can be written in shorthand as  
	\begin{equation}
		F(\phi_N) = 0. \label{system-N}
	\end{equation}
This system can be solved by a standard iterative method, such as Newton's method. 
%
%
	In this system, the value of phase speed $c$ has to be fixed for computing one particular solution. 
	Such an approach becomes impractical when a turning point on the bifurcation curve appears. 
	
\begin{figure}
\centering
		\includegraphics[width=0.35\textwidth]{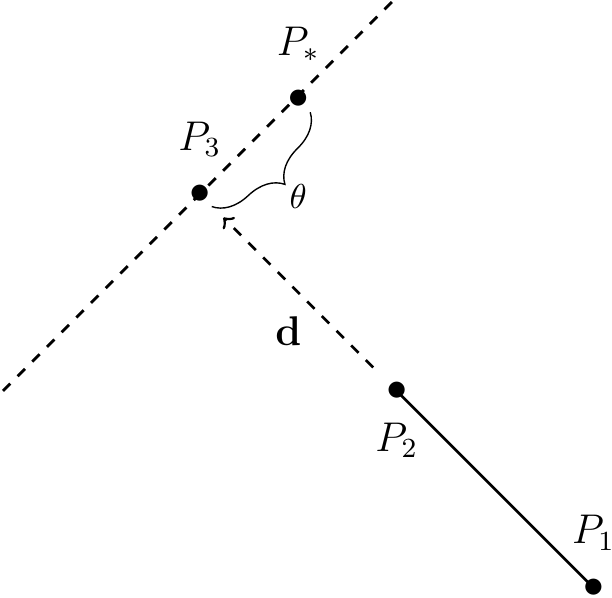}
\caption{\small Navigation on the bifurcation curve.}
\label{navig}
\end{figure}

	In \textsf{SpectraVVave} a different approach is employed: both the amplitude $a$ and the phase speed $c$ of a solution are treated as functions of a parameter $\theta$: $ a = a(\theta)$, $c = c(\theta)$. The parameter $\theta$ is to be computed from the system~\eqref{system-N3}. This construction makes is possible to follow turning points
on the bifurcation branch with relative ease.
	Having computed two solutions, i.e., two points on the bifurcation curve $P_1 = (c_1, a_1)$ and $P_2 = (c_2, a_2)$, one may find a direction vector $\mathbf{d}=(d^c, d^a)$ of the line that contains these points:
	\begin{align}
		\mathbf{d}:	\quad d^c = c_2 - c_1, \qquad	d^a = a_2-a_1.
	\end{align}
	Then the point $P_3 = (c_3,a_3)$ is fixed at some (small) distance $s$ from the point $P_2$ in the direction $\mathbf{d}$.
	\begin{align}
		P_3: \quad c_3 = c_2 + s \cdot d^c, \qquad	a_3 = a_2 + s \cdot d^a.
	\end{align}
	 The point $P_3$ plays the role of the initial guess for velocity and amplitude when computing the next solution $P_\ast = (c_\ast, a_\ast)$. 
	The solution point $P_\ast$ is required to lay on the line with direction vector $\mathbf{d_{\bot}}=(d^c_{\bot}, d^a_{\bot})$, which is orthogonal to the vector $\mathbf{d}$.
	\begin{align}
	\mathbf{d_{\bot}}: \quad d^c_{\bot} = - d^a, \qquad &d^a_{\bot} = d^c,  \\
	P_\ast: \quad c_\ast = c_3 + \theta d^c_{\bot} \qquad &a_\ast = a_3 + \theta d^a_{\bot}.
	\end{align}
	 A schematic sketch of finding a new solution $P_\ast$ is given in \autoref{navig}.	
	
	The variable $\theta$ is computed by Newton's method from the extended system 
\begin{equation}
F\left(
\begin{matrix}
\phi_N(x_1)\\ \vdots \\ \phi_N(x_N) \\ B \\ \theta
\end{matrix}
 \right)= \left(
\begin{matrix}
(-c+\mathcal{L}_N)\phi_N(x_1)+ f(\phi_N(x_1)) - B
\\ \vdots 
\\ (-c+\mathcal{L}_N)\phi_N(x_N)+ f(\phi_N(x_N)) - B
\\ \Omega(\phi_N, c, a, B)
\\ \phi_N(x_1)-\phi_N(x_N)-a
\end{matrix}
 \right)=
\left(
\begin{matrix}
0
\\ \vdots 
\\ 0
\\ 0
\\ 0
\end{matrix}
 \right). \label{system-N3} 
\end{equation}
		
	Here, a nonhomogeneous problem ($B \neq 0$) is considered.
	The equation 	
\begin{equation}
	\phi_N(x_1)-\phi_N(x_N)-a = 0, \label{amplitude_def}
\end{equation}
	 makes the waveheight of the computed solution to be that of $a$. The equation 
\begin{equation}
	\Omega(\phi_N, c, a, B)=0, \label{boundary}
\end{equation}
	 is called the  \emph{boundary condition}. 
	It allows to enforce different specifications on the computed traveling wave solution.   	
	For example, if one sets  
\begin{equation}
	\Omega(\phi_N, c, a, B)= \phi_N(x_1)+\cdots+\phi_N(x_N), \label{average}
\end{equation}
	 then \emph{the mean of the computed wave over a period will have to be equal to zero}.
	One may also experiment with	
\begin{equation}
	\Omega(\phi_N, c, a, B) = B, \label{homogeneous}
\end{equation}
	 to consider the homogeneous problem $(B=0)$.
	It can be also interesting to set
\begin{equation}
	\Omega(\phi_N, c, a, B) = \phi_N(x_N). \label{solitary}
\end{equation}
	 This enables us to compute traveling wave solutions that mimic solitary wave solutions. 	

\subsection{Convergence.}
In order to test the numerical implementation of the discretization,
the method is applied to a case where the solution is known.
One such case is the KdV equation
\begin{align*}
u_t + u_x + \frac{3}{2} u u_x + \frac{1}{6}u_{xxx} = 0, 
\end{align*}
which has a known solution, given in the form
\begin{align}
u_{\mathrm{exact}}(x,t) = a \sech^2\Big( {\textstyle\sqrt{\frac{3a}{4}}} (x-c t)\Big),\\
\label{KdV_sol}
\end{align}
with $c=1+a/2$.
%
Using the boundary equation~\eqref{solitary}, \textsf{SpectraVVave} is capable of computing 
approximations to solitary wave solutions of nonlinear wave equations. 
Solitary wave solutions are treated as traveling waves with sufficiently long wavelength 
that have the wave trough at zero.
In case of the KdV equation solitary wave solutions have exponential decay, 
and therefore, considering the symmetry of solitary solutions, 
the half-wavelength of $30$ is considered for the comparison. 
Approximation errors are summarized in \autoref{t2}.
%
\begin{table}	[ht]
\centering
\begin{tabular}{rrrr}
  \toprule
Nb.\ of grid points & $\log_{10}(\|u_{\mathrm{exact}} - u\|_{L^{\infty}})$ & $\log_{10}(\|u_{\mathrm{exact}} - u\|_{L^2})$ & Ratio of $L^2$-errors\\
\midrule
32 & $-1.389$  & $-2.092$ & \\
  
64 & $-3.705$  & $-4.549$ & $286.8$\\
  
128 & $-8.809$  & $-9.508$ & $90935.0$\\
  
256 & $-15.353$  & $-16.144$ & $4329670.9$\\
  
512 & $-15.353$  & $-16.087$ & $0.9$\\
\bottomrule
\end{tabular}
\caption{\small Estimates of error between the exact and computed solitary wave solutions for the KdV equation. Half-wavelength $l=30$, waveheight $a=1.2651$.}
\label{t2}
\end{table}

\begin{figure}[ht]
\centering
\includegraphics[width=0.47\textwidth]{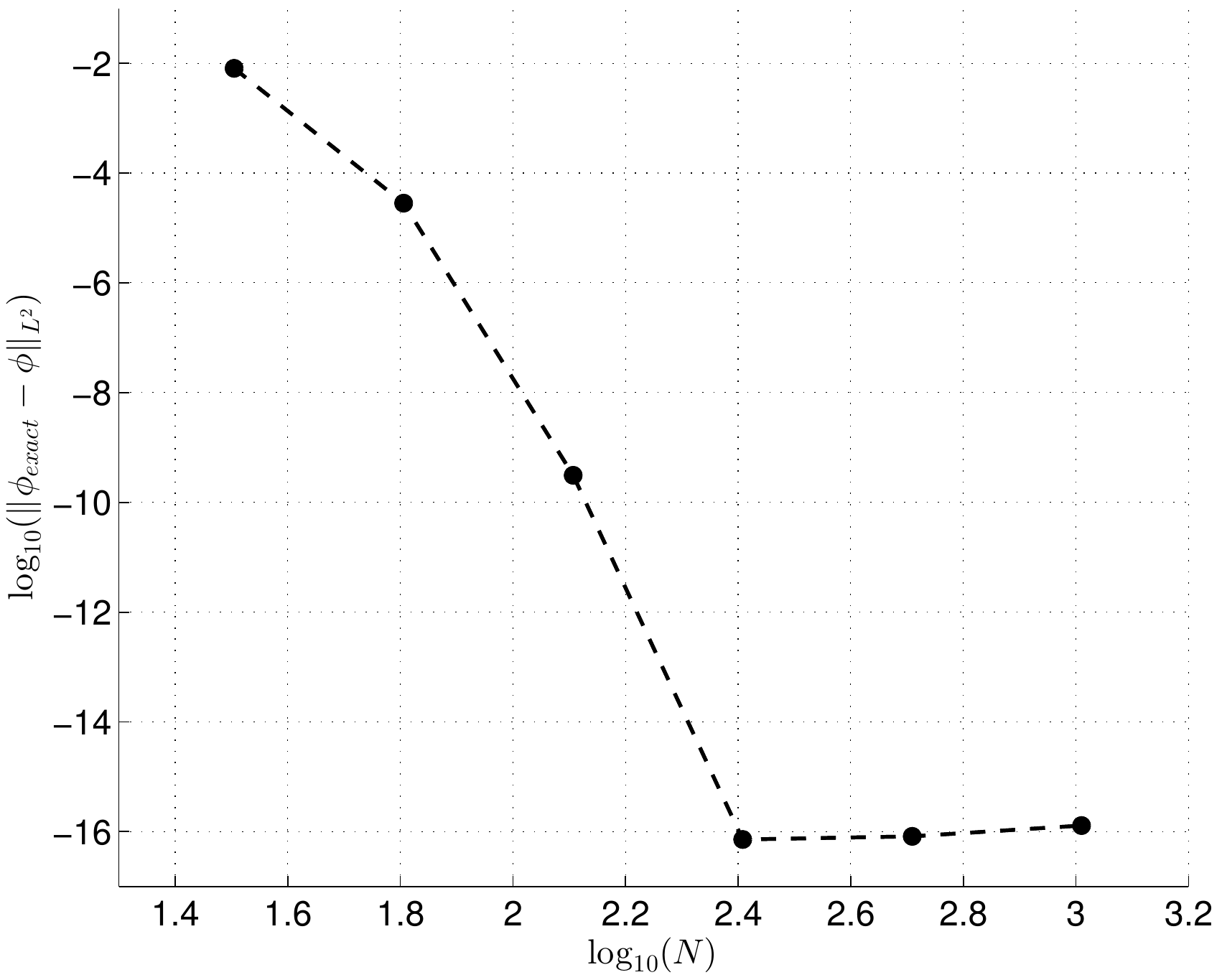}      
\caption{\small Graph of error estimates given in \autoref{t2}.}
\label{L2convergence}
\end{figure} 

\section{Experiments with SpecTraVVave.}
\label{sec:experiments}

\subsection{Termination of the waveheight-velocity bifurcation curve of the Whitham equation.}
\label{sec-W-term}
The waveheight-velocity bifurcation curve of the Whitham equation 
\begin{align}
	u_t + \frac{3}{2} u u_x + Ku_x = 0, \qquad 
	\widehat{Ku}(k) = \sqrt{\frac{\tanh(k)}{k}},  \label{W}
\end{align}
was studied numerically in~\cite{EK2}. 	
An attempt was made to identify the termination point of the Whitham bifurcation curve. 
The investigation was limited by computational tools and complete results were not obtained. 
In particular, the authors could not confirm that traveling wave solutions do not exist 
past the point where the authors, based on pioneering work of Whitham
\cite{Wh1} suspected a cusped solution.
In this section a number of numerical results on nature of the bifurcation curve
for the Whitham equation are presented.
Solutions to the equation~\eqref{W}
are computed in the form of traveling waves $u(x,t) = \phi(x - c t)$
and the homogeneous $(B = 0)$ integrated version the equation is considered:
\begin{align}
-c \phi + \frac{3}{4} \phi^2 + K \phi = 0. \label{W-int}
\end{align}
Special attention is given to relation between stability of solutions 
and their waveheight and velocity parameters, i.e., their position on the bifurcation curve. 
The following questions are under study: 
\begin{enumerate}[label=\alph*)] 
\item Where does the bifurcation curve terminate?

\item Where on the bifurcation curve do solutions change their stability?

\item Is there any role that the turning point on the bifurcation curve plays?
\end{enumerate}

\begin{figure}[ht]
\centering
\includegraphics[width=0.48\textwidth]{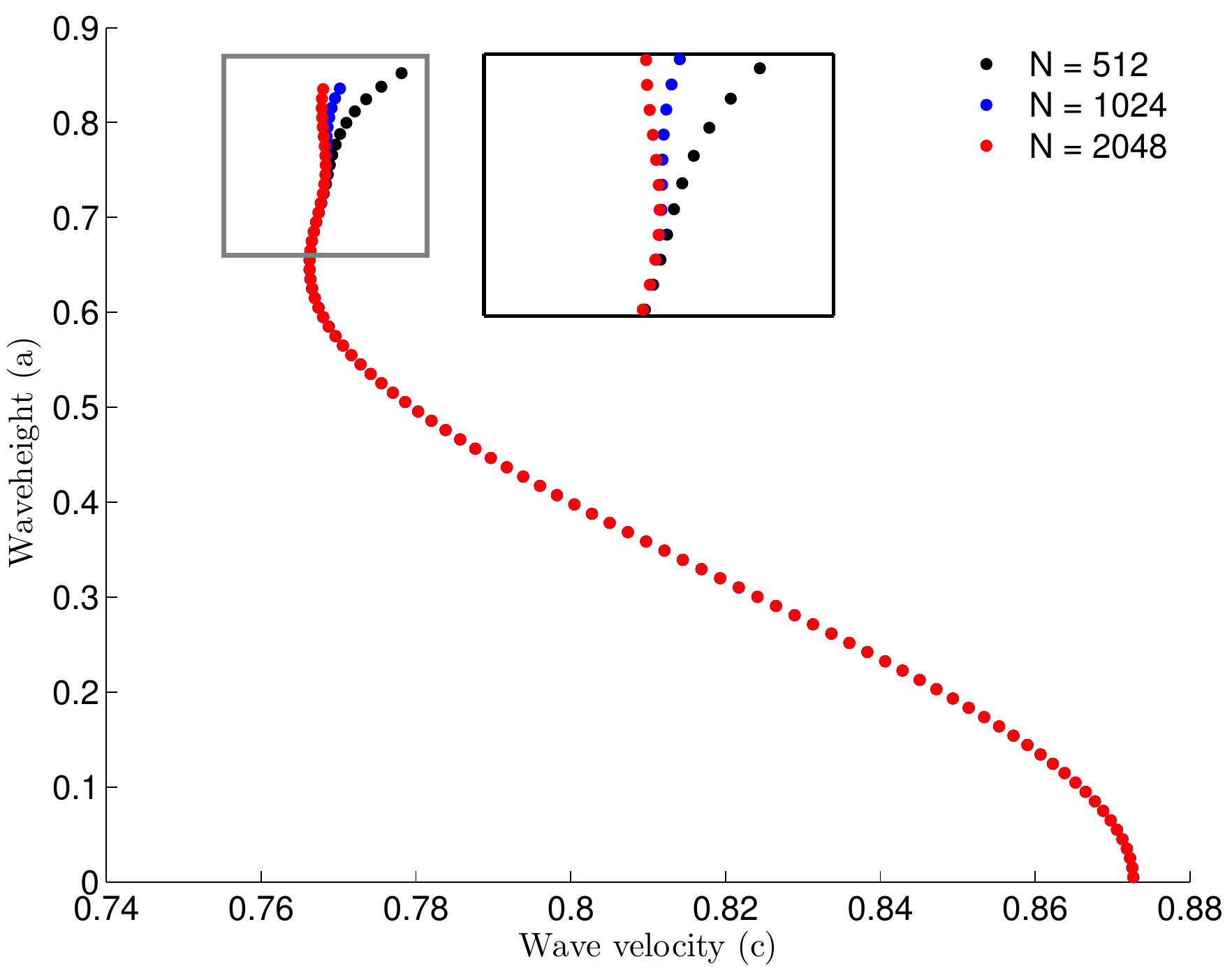}{}      
\caption{\small Whitham bifurcation curve in different grid resolutions.}
\label{bif-curve-all}
\end{figure}

The results presented here focus on $2\pi$-periodic solutions to the equation~\eqref{W-int}, 
i.e., solutions of the system~\eqref{system-N3}. 
\autoref{bif-curve-all} presents Whitham bifurcation curves with numbers of grid points 
$N = 512$, $N=1024$ and $N=2048$. 
The current implementation of the \textsf{SpecTraVVave} package allows fixing the number of grid points $N$ 
and a so-called doubling parameter $\DD$, 
i.e., the number by which $N$ is doubled as computations are made. 	
This allows us to get sets of solutions with $N, 2N, \ldots, 2^{\DD} N$ grid points. 
If $\DD = 1$ then only two sets of solutions are computed and they are regarded as 
lower grid (lower resolution) and higher grid (higher resolution) solutions.
While the system~\eqref{system-N3} is processed by Newton solver, 
lower grid solutions are taken as initial guesses for higher grid solutions. 
All curves shown in this manuscript have been produced
after tests with a number of resolutions were run, 
and the curves shown did not change significantly under further refinement.


\begin{figure}[ht]
\centering
\subfigure[][Bifurcation curve in the wavespeed-waveheight parameter space]{\includegraphics[width=0.48\textwidth]{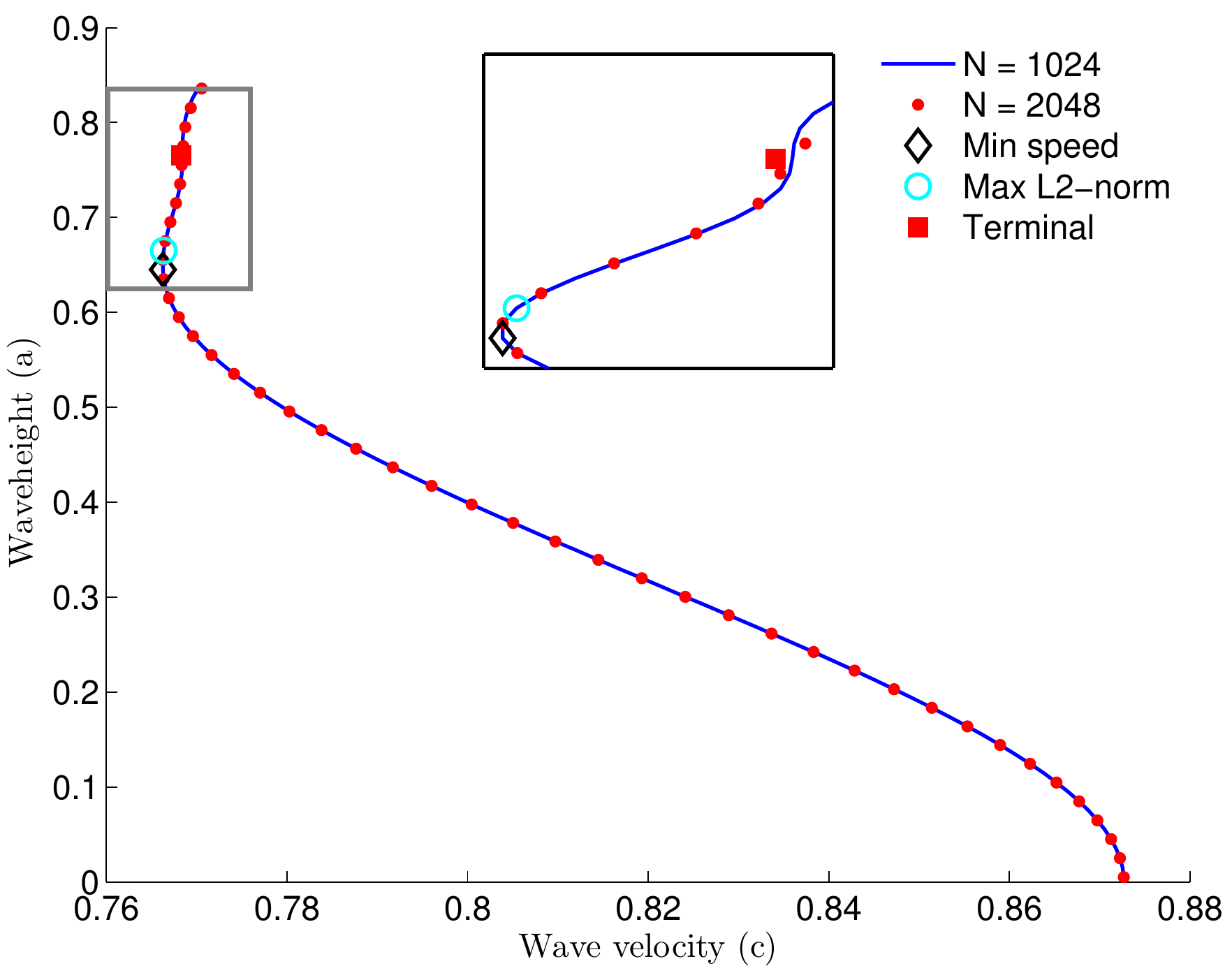}{}\label{sfig:bif-ch}}%
~%
\subfigure[][Bifurcation curve in a wavespeed-$L^2$-norm diagram]{\includegraphics[width=0.48\textwidth]{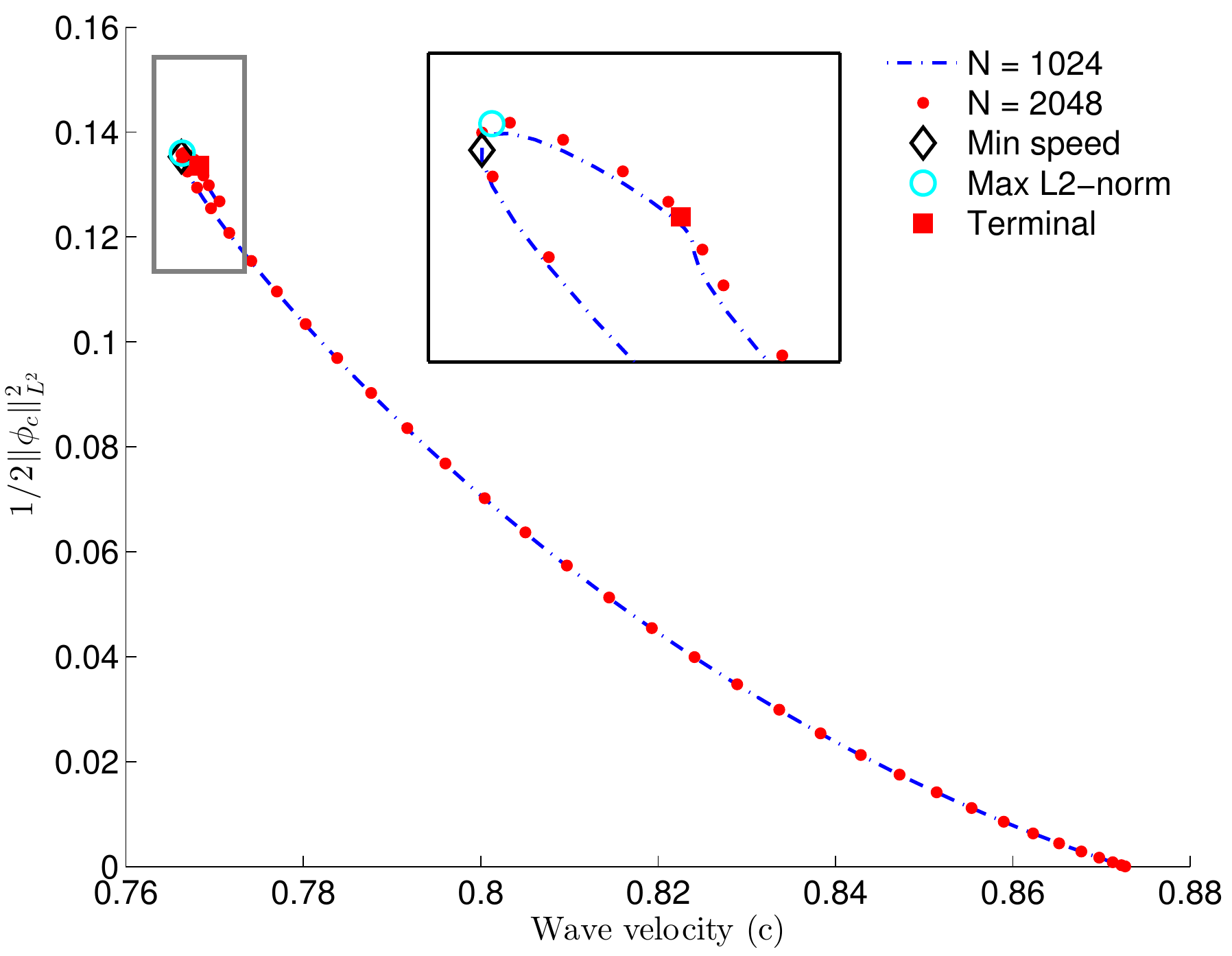}{}} 
\caption{\small
Whitham bifurcation curves for $2\pi$-periodic solutions. 
}
\label{bif-L2-1024}
\end{figure}

\begin{figure}[t]
\centering
\subfigure[][Profiles of the solutions singled out in \autoref{bif-L2-1024}\subref{sfig:bif-ch}]{\includegraphics[width=0.48\textwidth]{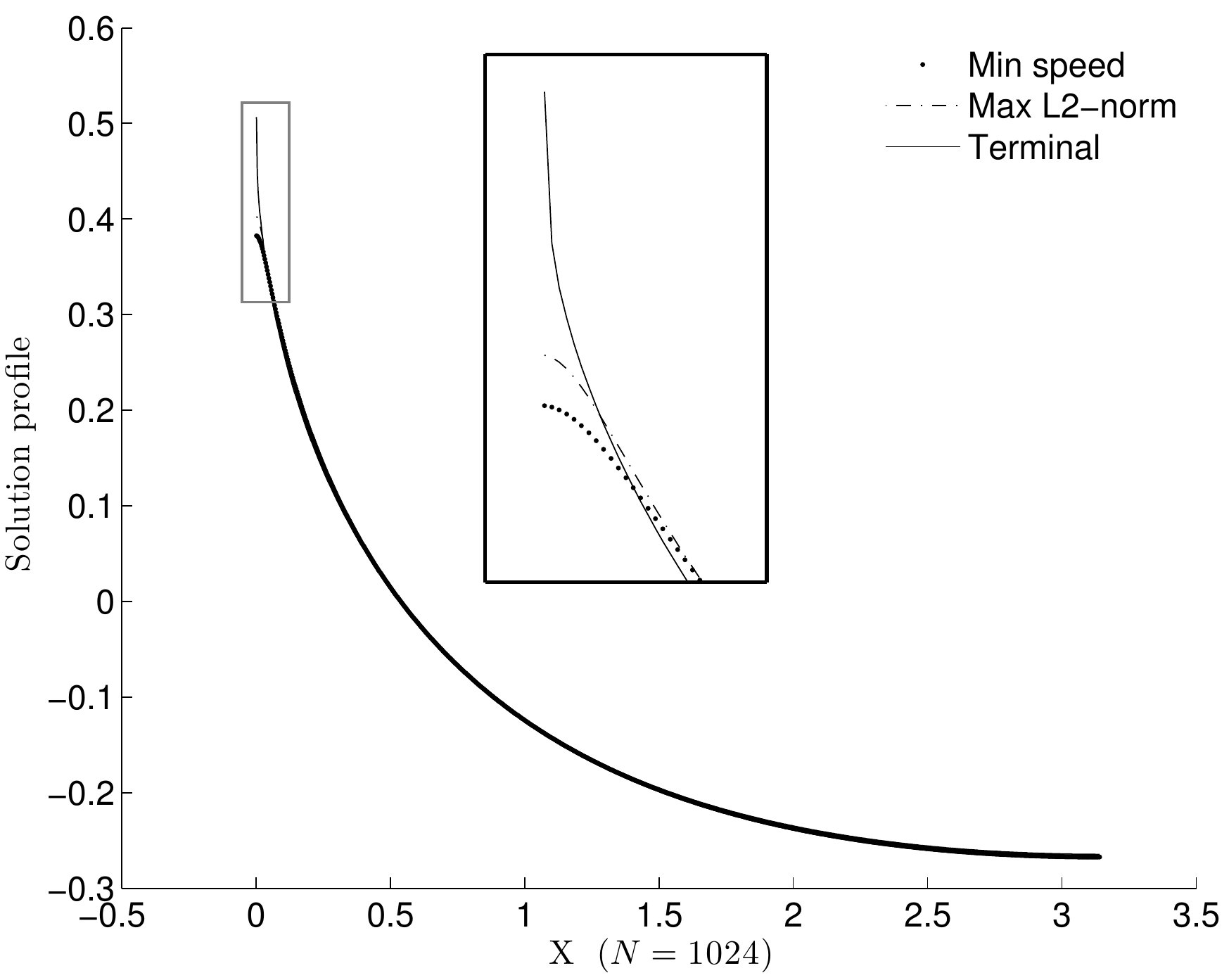}{}\label{sfig:prof-sing}}
~%
\subfigure[][Profile after terminal solution]{\includegraphics[width=0.48\textwidth]{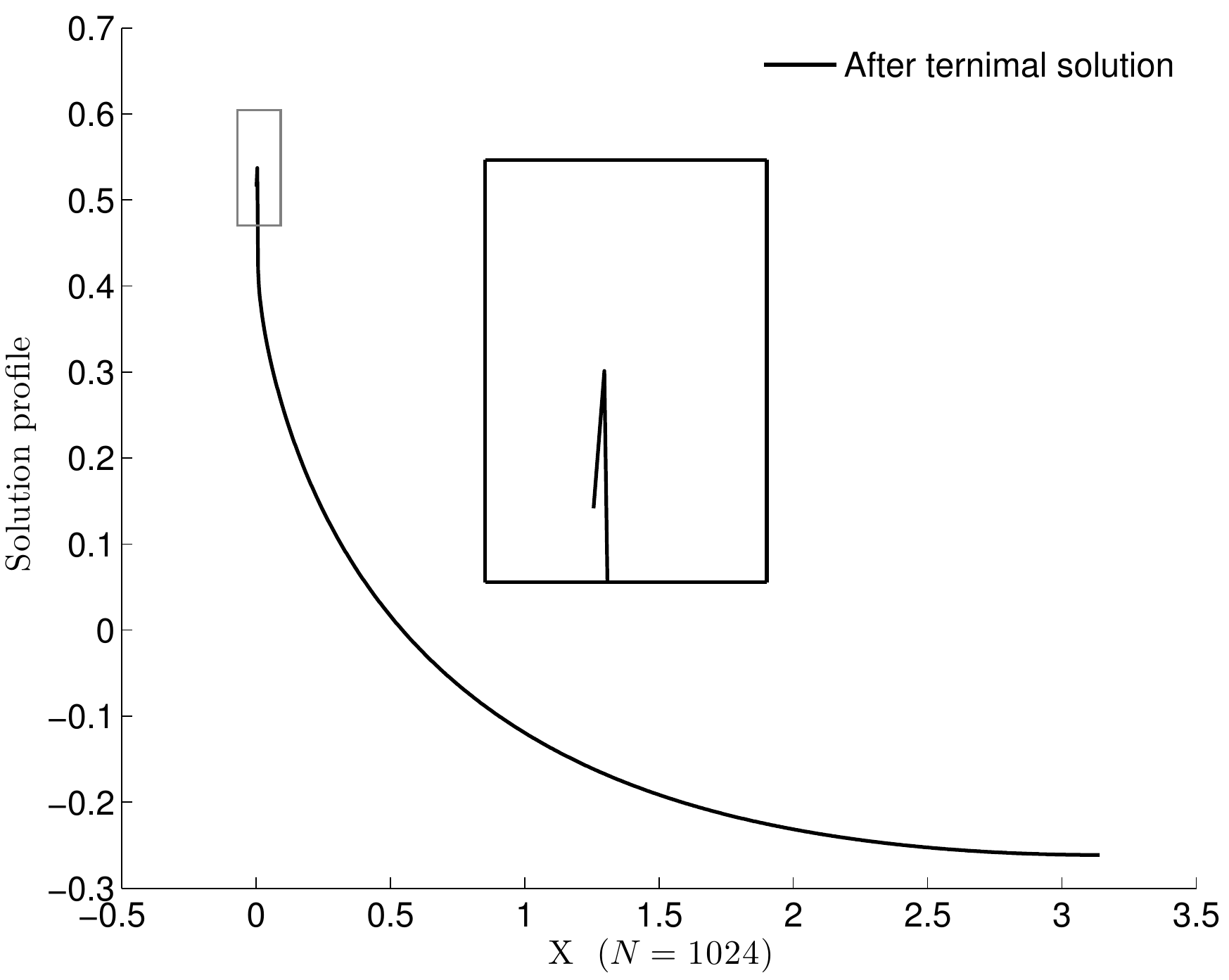}{}\label{sfig:prof-term}}
\caption{\small Profiles of specific traveling wave solutions.}
\label{3waves}
\end{figure}

\autoref{bif-L2-1024}\subref{sfig:bif-ch} 
presents the Whitham bifurcation curve computed by \textsf{SpectraVVave} 
with $N=1024$ and $\DD=1$. 
There are three solutions which deserve to be singled out:
\begin{enumerate}
	\item Traveling wave solution with minimum velocity (rhombus);
	\item Traveling wave solution with maximum $L^2$-norm (circle);
	\item Cusped traveling wave solution (square).
\end{enumerate}
Profiles of the above listed solutions are given in \autoref{3waves}\subref{sfig:prof-sing}.
The solution with minimum velocity corresponds to the turning point of the bifurcation curve.
The solution with maximum $L^2$-norm is very close to the latter one, although it has a higher waveheight 
and a different velocity.
The solution marked by a square is called here the \emph{terminal solution}.
As already mentioned, previous studies, such as~\cite{EK1,EK2} did not provide any conclusive
analysis on the part of the bifurcation curve past the turning point.
In particular, it was not clear whether solutions ceased to exist
at or after the turning point, or whether solutions were stable
or unstable after the turning point. 

Let us first focus on the stability of solutions.
Note that \textsf{SpectraVVave} has an evolution integrator routine, 
which enables one to check the stability of computed solutions. 	
The current version of the package uses the fourth-order method developed in~\cite{dFSS}.
In addition one may use a more refined analysis, resting on the evaluation
of invariant functionals. This analysis is based on the the observation
that the traveling waves can be thought of as solutions of a constrained
minimization problem. This analysis is based on ideas developed by
Boussinesq, first exploited in~\cite{B2}, and later used 
in~\cite{BSS,TJB,NgK}, and many other works.
Let us define two functionals $V$ and $E$:
\begin{align}
	V( \phi) = \frac{1}{2} \int_{-\infty}^{+\infty} \phi^2\,\dd\zeta, \qquad
	E( \phi) = \int_{-\infty}^{+\infty} \left\{ \Sfrac{1}{2} \phi^3 - \phi \, K\!  \phi \right\}\, \dd\zeta.
\end{align}
The equation~\eqref{W-int} can be then written in terms of variational derivatives of $E$ and $V$
as
\begin{align}
	E' ( \phi ) - c V ' ( \phi ) = 0.
\label{zero}
\end{align} 
It is known from~\cite{BSS} that the stability of solitary wave solutions depends 
on convexity of the function $d(c) = E(\phi) - c V(\phi)$. 
Solutions with values of $c$ for which $d''(c) > 0$  are stable solutions,
and solutions with wave speeds for which $d''(c) <0$ are unstable solutions.
	
If the current numerical investigation confirms the latter hypothesis, 
then the equation~\eqref{d'c} establishes a direct relation 
between the stability of traveling wave solutions and their $L^2$-norms, 
in case of the Whitham equation.
Note that differentiation of $d(c)$ yields
\begin{align}
	d '(c) &= \underbrace{E'(\phi) - c V'(\phi)}_{=0} + V(\phi).
\end{align}
Using~\eqref{zero} as indicated yields
\begin{align}
	d '(c) &= V(\phi) =  \frac{1}{2} \int_{-\infty}^{+\infty} \phi ^2\,\dd\zeta = \frac{1}{2}\|\phi\|_{L^2}^2 . 
\label{d'c}
\end{align}
Therefore, in order to understand the convexity of $d(c)$, it is sufficient
to find points of maximum $L^2$-norm on the curve in the right panel
of \autoref{bif-L2-1024}.
It is straightforward to see that $d''(c)$ changes sign in the neighbourhood of the maximum point 
of this curve, i.e., around the solution with maximum $L^2$-norm. 
In particular, $d''(c) > 0$, i.e., solutions are stable to the left of the maximum point, 
and $d''(c) < 0$, i.e., solutions are unstable to the right of the maximum point.
	
\begin{figure}[t]
  \subfigure[][Minimum speed solution]{\includegraphics[width=0.31\textwidth]{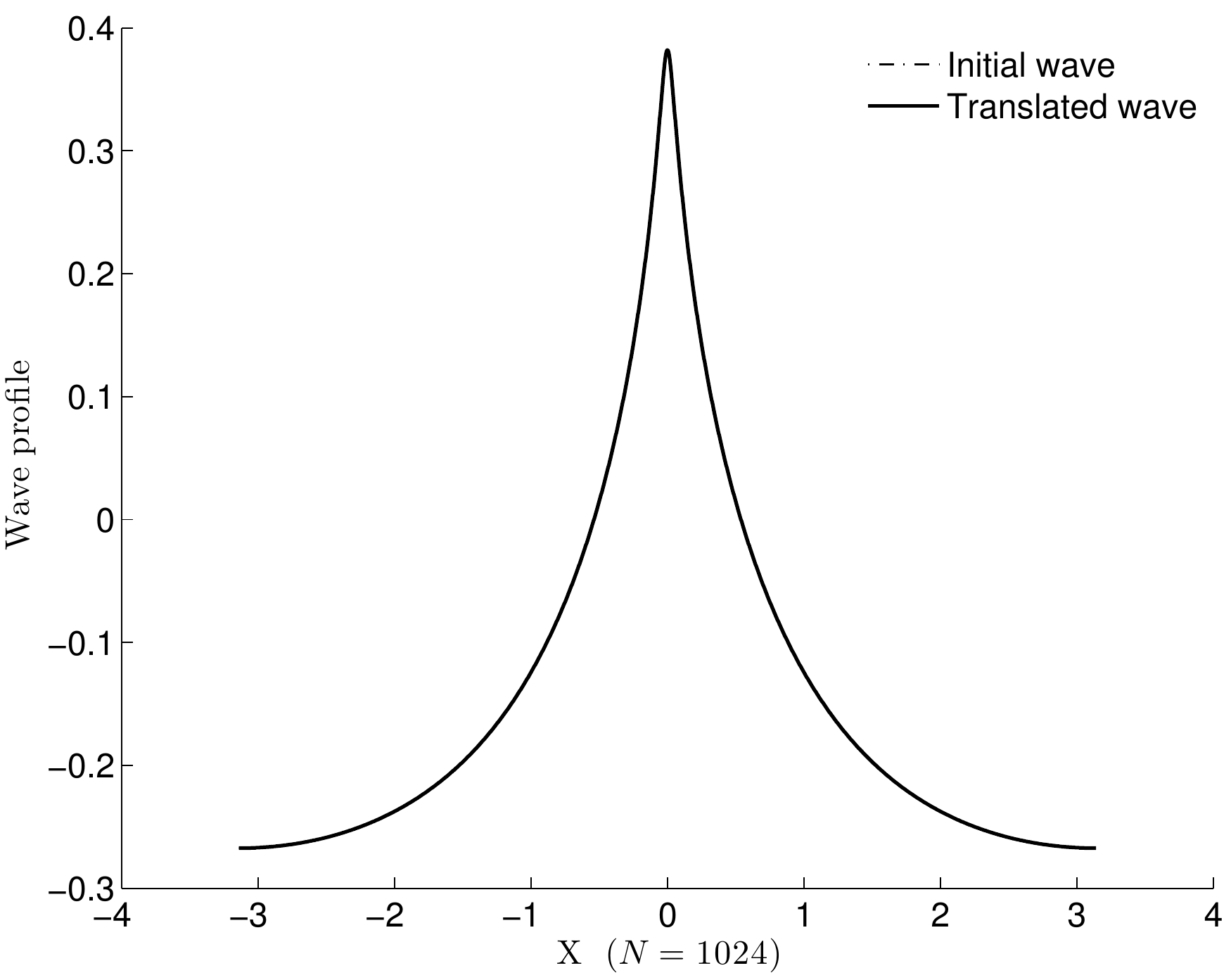}}
~%
\subfigure[][Maximum $L^2$-norm solution]{\includegraphics[width=0.31\textwidth]{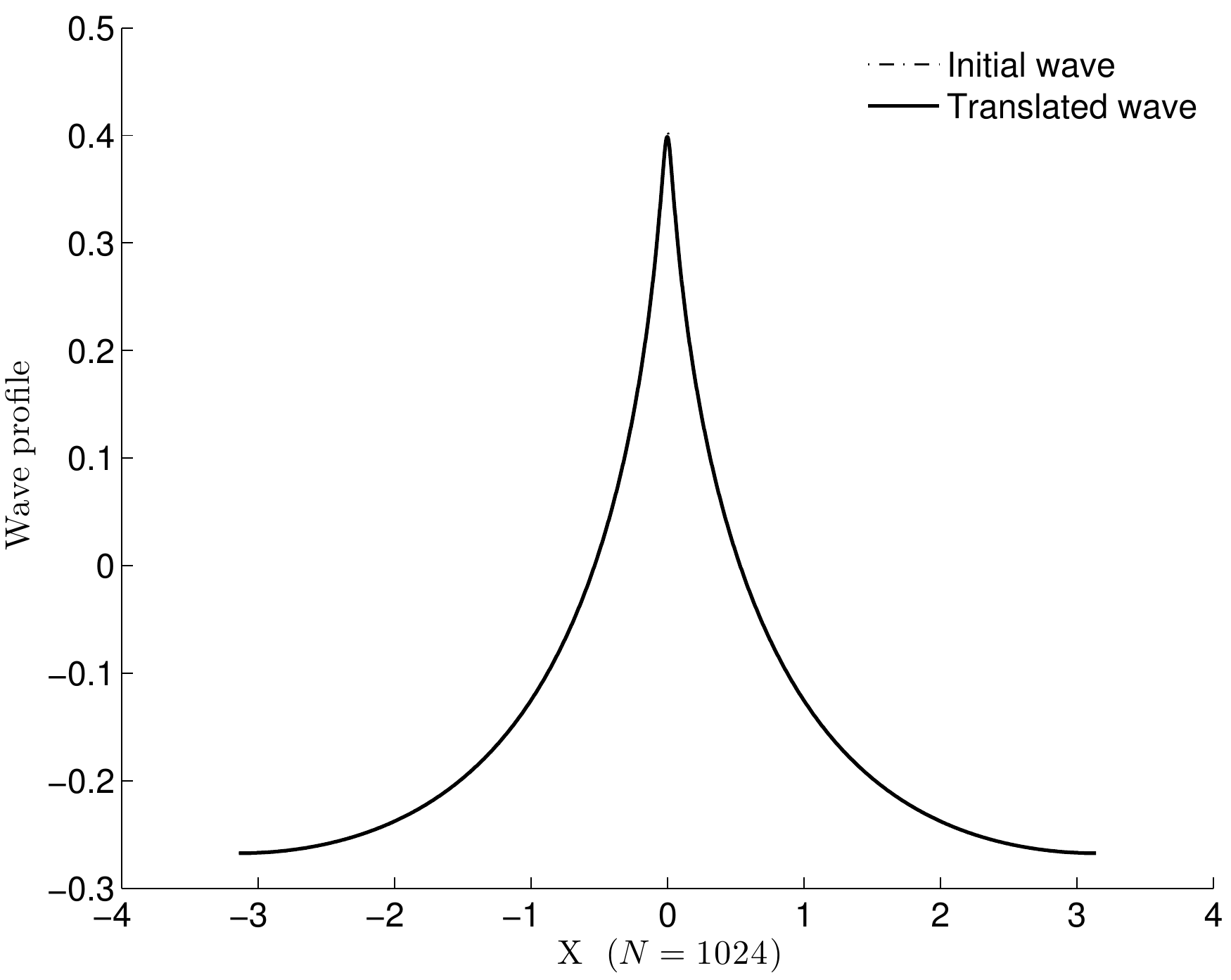}{}}
~%
\subfigure[][Terminal solution]{\includegraphics[width=0.31\textwidth]{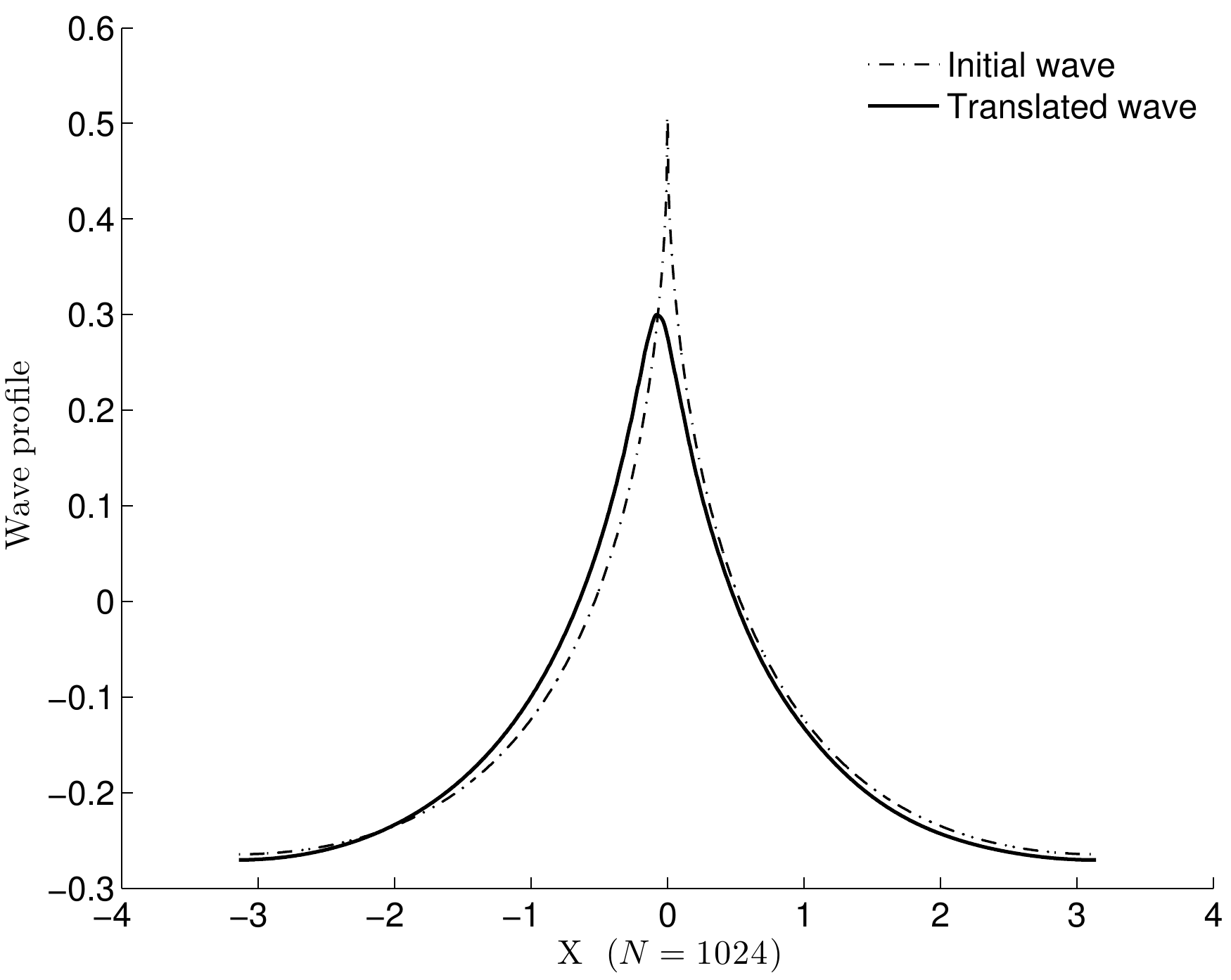}{}}
\caption{\small Evolution of specific solutions in time. Time range for each solution is three periods.}
\label{time-evo}
\end{figure}

In addition, the solutions were tested with the evolution integrator to confirm their stability/instability in time. 
The solution with maximum $L^2$-norm and those on the left-hand side of it seemed to be stable in time. 
Solutions on the right-hand side do not preserve their shape and thus are unstable. 
Examples are given in \autoref{time-evo}. 
This analysis confirms that 
the point corresponding to the minumum wave speed (the turning point), 
and the point of stability inversion are two distinct points on the bifurcation curve.
Moreover, the point of stability inversion
is a little further up the branch from the turning point.

Next, we turn our attention to the analysis of the terminal point. There
are two main questions. Does the branch terminate, and if so, 
does the terminal point on the branch correspond to a cusped traveling wave.
First of all, note that the solution, which is computed by \textsf{SpectraVVave}, 
past the terminal solution has two crests, no matter how small the 
stepping on the bifurcation branch is taken.
(see \autoref{3waves}\subref{sfig:prof-term}).   
Secondly, as will be explained presently, the relation 
\begin{equation}
\frac{c}{\sup_{x} \phi(x)} = \frac{3}{2} 
\end{equation}
holds for the terminal solution with a good degree of approximation.
For the most accurate runs, we obtain $c / \sup_{x \in \R} \phi(x)$ $\approx$ $1.51$. 
To explain how this relation comes about note that the steady integrated form of the
Whitham equation can be written as
\begin{equation}\label{Whitham_Alt}
\Big( \frac{c}{\sqrt{3}} - \frac{\sqrt{3}}{2}\phi \Big)^{2}  = \frac{1}{3} c^2 - K \phi.
\end{equation}
It is clear that for any $\phi < \frac{3}{2} c$, the relation~\eqref{Whitham_Alt}
can used in a bootstrap argument to show that any continuous solution must be
in fact smooth. However for the case $\phi =  \frac{3}{2} c$ this bootstrap
argument fails since the left-hand side vanishes. It can be concluded that
a solutions containing a cusp will have a maximum  value of $\frac{3}{2} c$.

	
As an additional check, the discrete cosine coefficients of the solutions were
examined, and fitted to the following models:
\begin{align}
	\mathcal{E}(k) = \nu_1 e^{-\nu_2 k^n}, \qquad 
	\mathcal{P}(k) = \frac{\mu_1}{\mu_2 + \mu_3 k^m},
\end{align}
%
where $\nu_1$, $\nu_2$, $\mu_1$, $\mu_2$, $\mu_3$, $n$ and $m$ are fitting parameters. 
A smooth function is known to have discrete cosine coefficients with exponential decay in $k$. 
On the other hand, if a function is not smooth, the discrete cosine coefficients feature polynomial decay. 
To identify the best fit, two parameters were used: $L^2$ norm of residual 
and Akaike information criterion (AIC) measure. 

From the data given in the \autoref{t1}, one can deduce that for solutions with minimum speed 
and maximum $L^2$-norm exponential fit is better than polynomial. 
That is not the case for the terminal solution. 
Thus, the first two solutions are smooth and the terminal solution is nonsmooth. 
In fact, the polynomial fit is better than exponential for solutions that 
are between the maximum $L^2$-norm solution and the terminal solution.

The numerical evidence brought forward supports the conclusion 
that the Whitham bifurcation branch terminates at the terminal point
indicated in \autoref{bif-curve-all}. Of course, as already mentioned,
this conclusion has now also been reached using tools of mathematical analysis
\cite{EW2016}.

\begin{table}	
\centering
\begin{tabular}{l|rr|rr|rr|}
\cline{2-7}
& \multicolumn{2}{ c     }{Min speed solution} & \multicolumn{2}{ c     }{Max $L^2$-norm solution}& \multicolumn{2}{ c |  }{Terminal solution} \\
\multicolumn{1}{l|} {Model} &$\mathcal{E}(k)$&$\mathcal{P}(k)$&$\mathcal{E}(k)$&$\mathcal{P}(k)$&$\mathcal{E}(k)$&$\mathcal{P}(k)$ \\ 
\midrule
\multicolumn{1}{l|} {Residual's $L^2$ norm} & $5\times 10^{-5}$ &$4\times 10^{-3}$ & $7 \times 10^{-5}$ & $3\times 10^{-3}$ &  $6\times 10^{-3}$& $6\times 10^{-4}$ \\
\multicolumn{1}{l|} {AIC} & -543 & -321 & -529 & -333 & -298 & -416 \\ 
\bottomrule
\end{tabular}
\caption{\small Results for measures of fit.}
\label{t1}
\end{table}
	

\subsection{Interaction of solitary wave solutions of modified Benjamin--Ono equation}
\label{mB-O-interact}

In this section, we utilize the \textsf{SpectraVVave} package to obtain
high-precision approximations to solitary-wave solutions of
the modified Benjamin--Ono
\begin{equation*}
u_t + u^2u_x + u_x - \mathcal{H} u_{xx} = 0,
\end{equation*}
which is a special case of the generalized Benjamin--Ono equation, with $p=2$.
This case corresponds to the critical scaling, i.e., invariance of 
the energy norm under the natural invariant scaling,
and was not investigated
in~\cite{HKJB} since it is more difficult than the supercritical cases,
where $p > 2$.

The Benjamin--Ono equation was found by Benjamin~\cite{B1} as a model
for long small-amplitude interfacial waves in deep water. 
The validity of approximating the more physically correct configuration of a continuous
density distribution by the two-layer approximation has recently been justified mathematically
\cite{ChenWalsh2016}.

Solitary-wave solutions of the modified Benjamin--Ono equation with $p=3$, $p=4$ and $p=5$
were approximated in~\cite{BoKa} with a standard Newton scheme.
The solutions in~\cite{BoKa} were not very accurate, but since singularity formation of the evolution
equations the accuracy of the solitary-wave approximation was not an important issue.
The problem with the method of~\cite{BoKa} and some other works was that the fft
used here was not purged of possible symmetries (translational and reflective).
In the current code, since a cosine formulation is chosen, these symmetries
are automatically eliminated, and the resulting computations are able to to
render more accurate approximations. 

Solitary-wave solutions of these equation could be computed with higher accuracy
using a type of Petviashvili method in~\cite{Pelinovsky}, but here we emply
the \textsf{SpectraVVave} package using the boundary equation~\eqref{solitary},
and treating solitary waves as 
traveling waves with sufficiently long wavelength 
that have wave trough at zero. 
Once these high-accuracy solutions are found, they are aligned in 
an evolution code using a high-order time integrator, and the
interaction of two waves is studied.

Two question are investigated. First, the interaction is investigated for evidence
of integrability. Second, we are looking for possible annihilation of one of the waves,
such as may happen in some other evolution equations~\cite{Tjon}.

\begin{figure}[t]
\centering
		\includegraphics[width=0.6\textwidth]{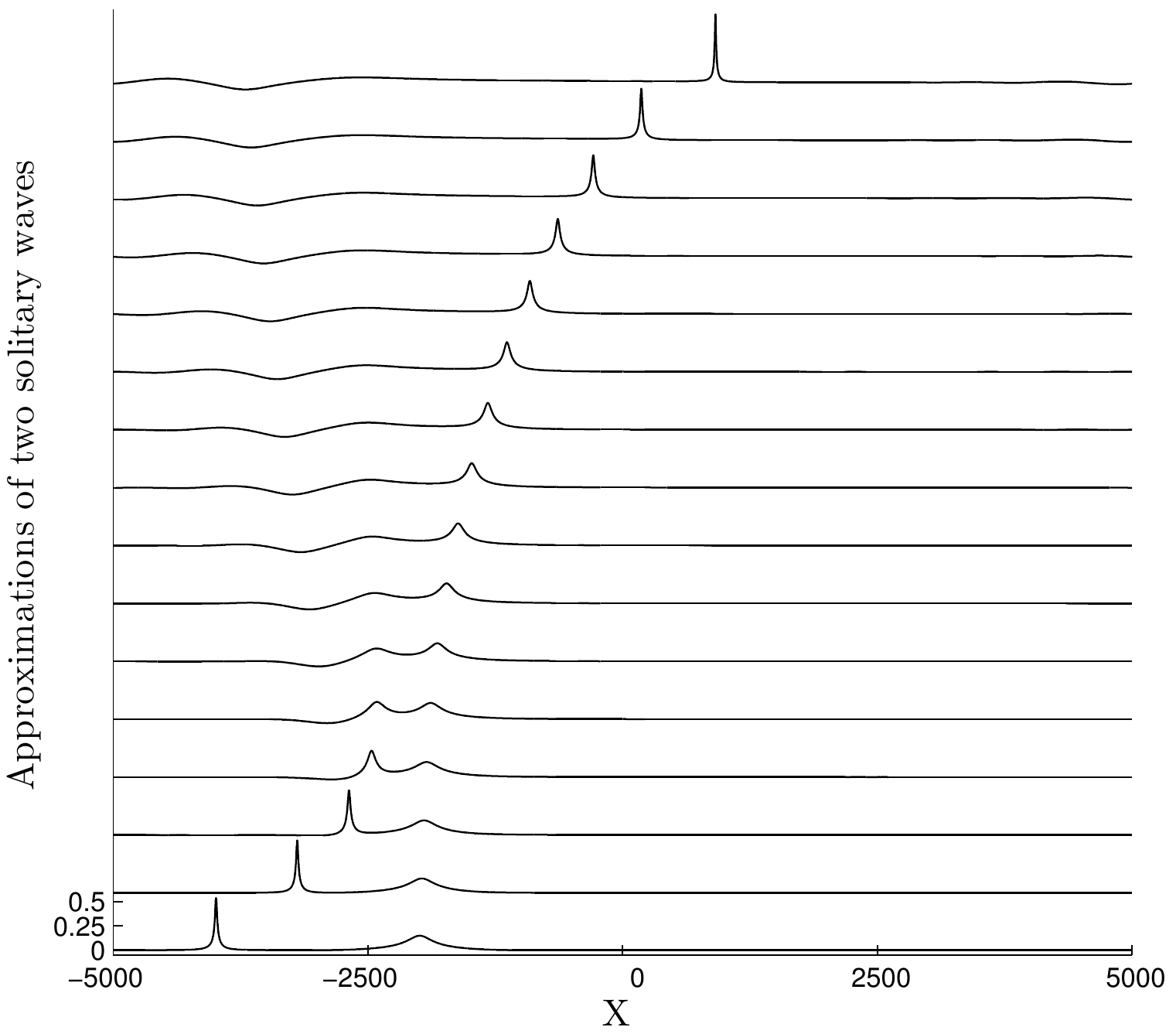}      
\caption{\small Interaction of two solitary waves of the modified Benjamin--Ono equation.}
\label{solitons-full}
\end{figure}

A possible approach to studying the question of complete integrability
is analyzing the interaction of two solitary wave solutions of the equation, 
such as carried out in~\cite{HK,HKJB} for other nonlocal equations.	
In the \autoref{solitons-full} snapshots of interaction of two solitary waves at different times are shown.
The time difference between two consecutive snapshots is constant.    
As it may be observed, during the process of interaction, the two initial solitary waves combine 
into a single wave, and an additional oscillation is produced. 
This leads us to the conclusion that the interaction of solitary waves is not elastic 
and the modified Benjamin--Ono equation may not be integrable. 
In addition, it appears that the smaller wave disappears as most of its
mass is acquired by the larger wave. Thus one may argue that the small
wave is annihilated by the larger wave.

\subsection{Effect of competing dispersion in the Benjamin equation}
\label{Benjamin_eq}
The Benjamin equation was found by Benjamin~\cite{B3} as a model for two-layer flow in the case
when the interface is subject to surface tension. The approximation may not be a good model for
a stratified situation, but more applicable to the case where two fluids are separated by a
sharp interface. 
The equation is 
\begin{equation}\label{benjamin}
u_t + u_x + uu_x - \mathcal{H} u_{xx} - \tau u_{xxx} = 0,
\end{equation}
where $\tau$ is a parameter similar to the Bond number in free surface flow~\cite{B3, KalischBenjamin, Walsh2014a}.
	

In this section, a study relating to the effects of competing dispersion operators 
on the shape of periodic traveling waves in the Benjamin equation is presented.
An in-depth study of solitary waves was carried out in~\cite{DDM}.
As will come to light, the periodic case features some new phenomena, such as secondary
bifurcations, connecting and crossing branches.
For the purpose of this study, we fix the parameter $\tau=0.1$,
so that the dispersion relation for the linearized equation is
\begin{align}
c(k) = 1.0 - |k| +  \frac{1}{10}  k^2. 	\label{Ben_disp_rel}
\end{align}

\begin{figure}[ht]
  \centering
\includegraphics[width=0.47\textwidth]{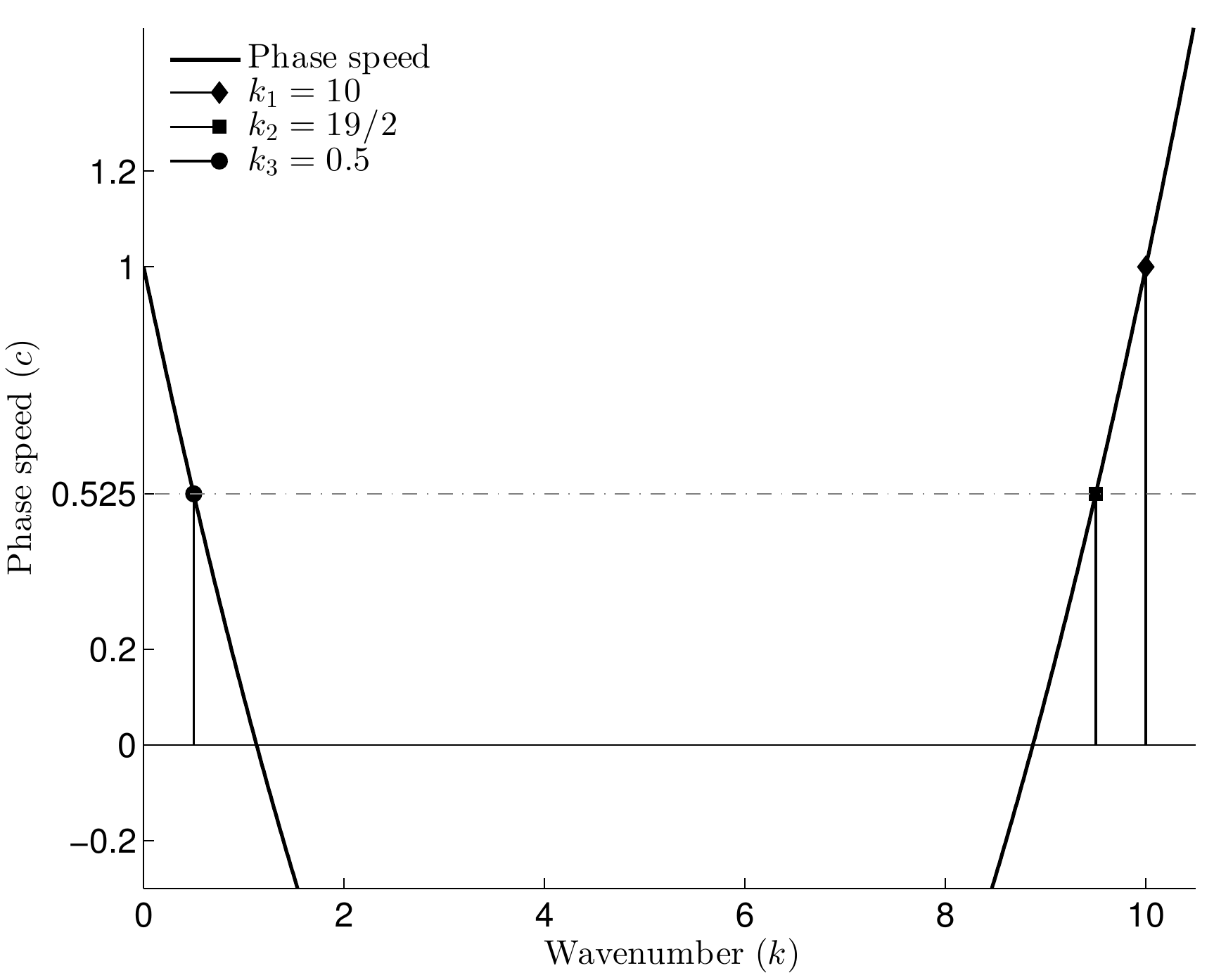}  \caption{Dispersion relation \eqref{Ben_disp_rel}}
\label{disp-rel}
\end{figure}

\begin{figure}
  \centering
\includegraphics[width=0.47\textwidth]{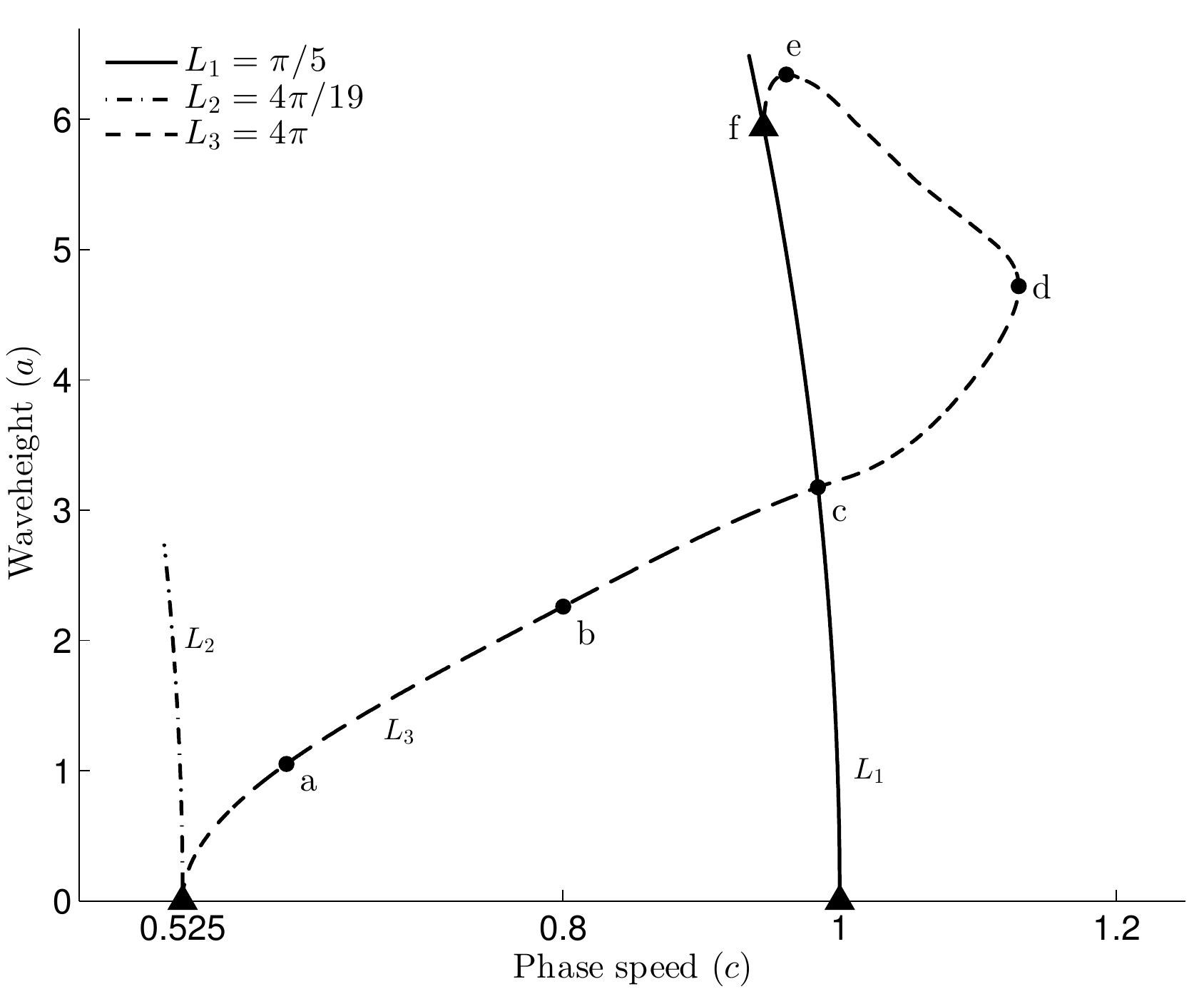}
\caption{Bifurcation curves of the equation~\eqref{benjamin} with different wavelengths, $\blacktriangle$ --  points of bifurcation, $\bullet$ -- selected solutions (see \autoref{long-solutions}).}
\label{fig:ben-bifur}
\end{figure}

Traveling wave solutions with fundamental wavelengths $L_1 = \pi/5$, $L_2 = 4\pi /19$ and $L_3 = 4\pi$ 
are computed for the equation~\eqref{benjamin}.  
The corresponding wavenumbers are $k_1 = 10$, $k_2 = 19/2$, and $k_3 = 0.5$, respectively.
A plot of the dispersion relation~\eqref{Ben_disp_rel} is given in \autoref{disp-rel}.
Bifurcation branches of traveling wave solutions with the selected wavelengths are given in \autoref{fig:ben-bifur}.

\begin{figure}[ht]
        \centering
		\includegraphics[width=0.47\textwidth]{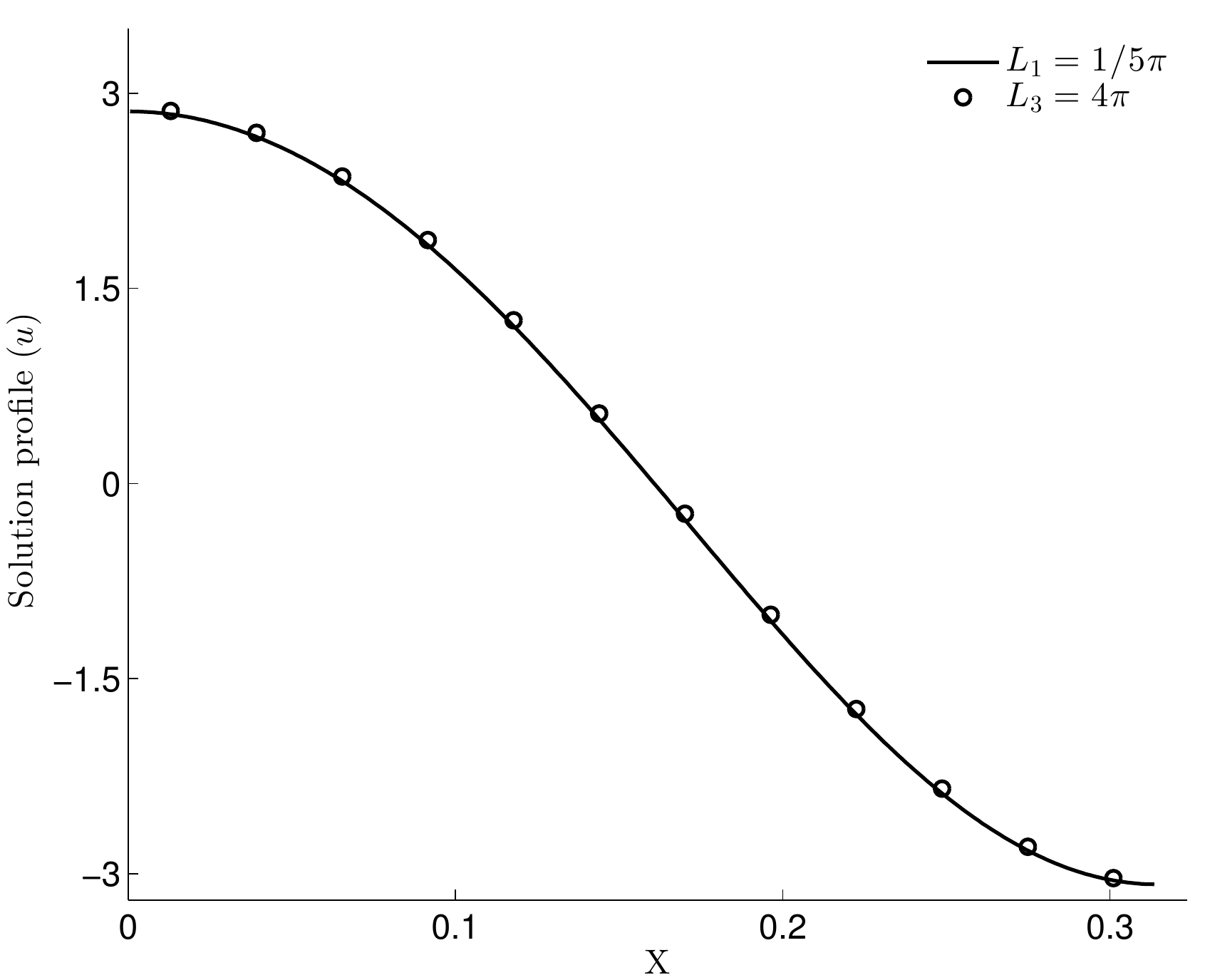}{}      
	\caption{\small Solution profiles at the point $(c^*, a^*)$ where the $L_1$ and $L_3$ branches meet (see \autoref{fig:ben-bifur}).}
\label{last-solutions}
\end{figure} 

\begin{figure}[ht]
\begin{subfigure}
	\centering
	\includegraphics[width=0.47\textwidth]{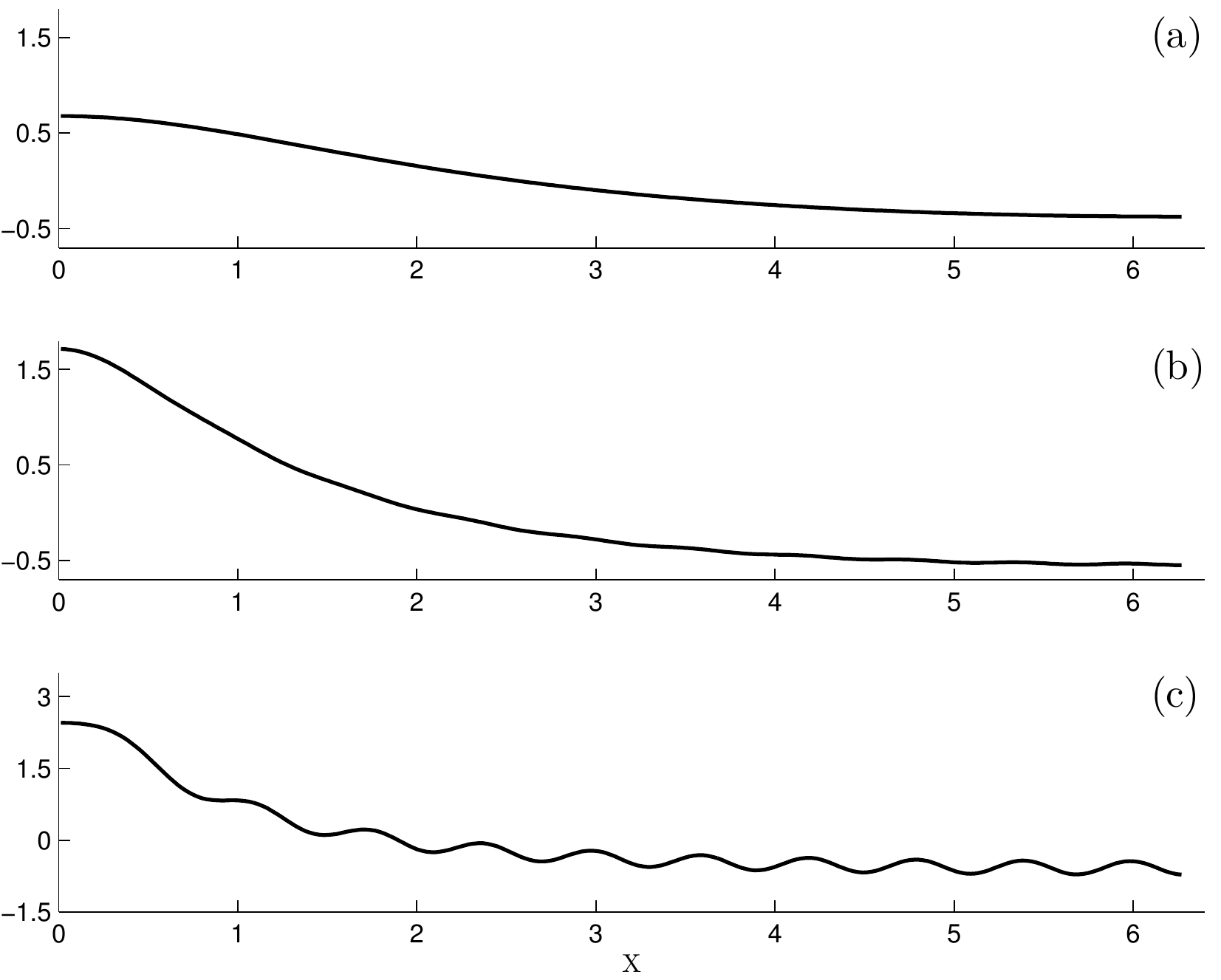}      
	    \end{subfigure}%
~%
\begin{subfigure}
        \centering
		\includegraphics[width=0.47\textwidth]{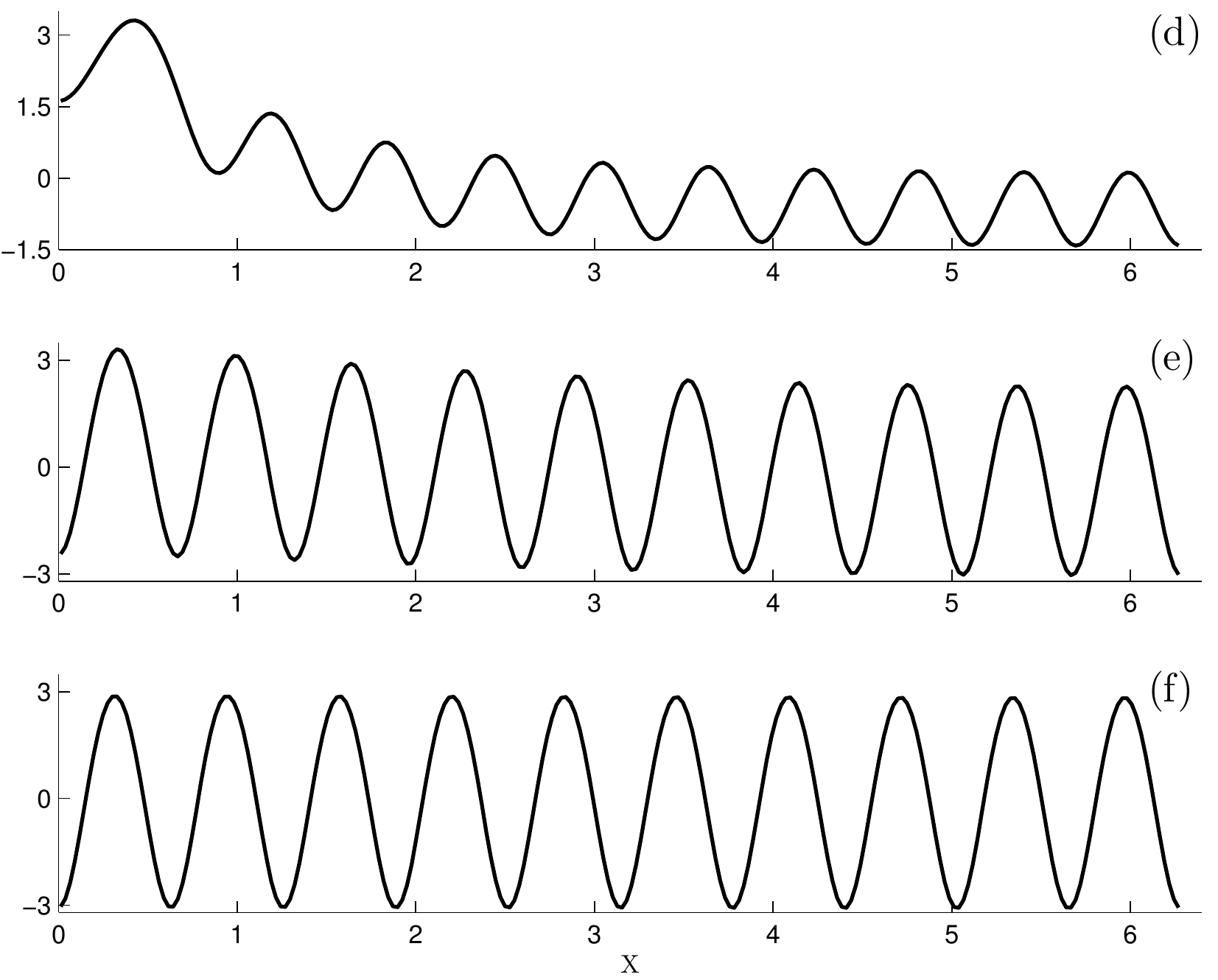}{}      
    \end{subfigure}%
	\caption{\small Selected solutions of the equation~\eqref{benjamin}. 
Labels preserved as shown in \autoref{fig:ben-bifur}.}
\label{long-solutions}
\end{figure}

The branch denoted by $L_1$ originates at the bifurcation point 
located at $c = 1$ and zero waveheight.
The branches denoted by $L_2$ and $L_3$ originate from the same bifurcation point,
located at $c = 0.525$ and zero waveheight. 
These two branches continue in different directions, due to differences in wavelength.
In particular, the $L_3$ branch contains waves with shorter wavelengths, and falls into the  
capillary regime. On the other hand, the $L_2$ branch falls in the gravity regime.
As the waveheight grows, solutions on the $L_3$ branch first cross the $L_1$ branch without
connecting. Additional oscillations develop in the solutions until 
a new fundamental wavelength $\pi/5$ is reached, and the branch terminates
as it connects to the $L_1$ branch. The situation is depicted in \autoref{long-solutions}. 
The point where the $L_1$ and $L_3$ branches meet
is approximately $(c^*,a^*) = (0.945, 5.938)$.
The corresponding profiles essentially overlap, as shown on \autoref{last-solutions}.
This point can also be interpreted as a secondary bifurcation point 
of the $L_1$ branch, where solutions with wavelengths that are multiples of $\pi/5$ develop. 
We should note that similar phenomena concerning crossing and connecting
branches were previously observed by Remonato and Kalisch~\cite{RK}
for the Whitham equation with surface tension which was introduced in~\cite{HurJohnson}. 
	
\section{Conclusions and future work}
\label{sec:conclusion}

The numerical algorithm of \textsf{SpectraVVave} features ample flexibility for researching 
different aspects of nonlocal dispersive wave equations and their traveling wave solutions.
The solver package is simpler in use when compared with programs such as \textsf{AUTO} and 
\textsf{Wavetrain}, however it does not have the same level of generality. 
Moreover, \textsf{AUTO} and \textsf{Wavetrain} are programmed in low-level programming languages
and will therefore run more efficiently.
\textsf{SpectraVVave} is implemented in an object-oriented fashion~\cite{PyBook}, 
which makes the program easily expandable.
\textsf{IPython} provides means for interactive work with the package, 
and enables users to create convenient notebook-programs.
A parametric approach in defining amplitude and phase speed makes it possible to follow
turning points on bifurcation curves. Specification of different boundary conditions allows 
computing solutions with certain features, such as traveling waves with mean zero, 
or approximations to solitary waves.

In this work, the \textsf{SpectraVVave} package has been put to use for the study of 
a number on nonliear evolution equations: the Whitham equation, the modified Benjamin--Ono equation
and the Benjamin equation.
For the chosen set of parameters, experiments on the Whitham equation resulted in numerical 
confirmation of the conjecture on cusped solutions. It was also possible to identify the point of stability 
inversion of traveling wave solutions of the equation and the termination point of its bifurcation curve.

In case of the modified Benjamin--Ono equation, the study on solitary wave solutions lead us to conclude that 
interaction process ended with annihilation of one of the two waves. 
The experiment on the Benjamin equation showed one more example of the effect of competing dispersion.
As the amplitude increased, traveling wave solutions of wavelength $4\pi$ developed
additional oscillations, and later connected up with a branch of solutions with
wavelength $\pi/5$.

Future work on the \textsf{SpectraVVave} package will be focused on development of its functionality and broadening the range of problems that can be studied. Possible extensions may include implementation of algorithms based on the Petviashvili method~\cite{Duran2014CNSNS}--\cite{Duran2014JCAM} and generalization to systems of equations.


\textbf{Acknowledgments.}
The authors would like to thank Mats Ehrnström and Erik Wahlén
for fruitful discussions on the subject of the current paper. 
This research was supported by the Research Council of Norway on grant no. 213474/F20.

\appendix
\section*{Appendices}
\addcontentsline{toc}{section}{Appendices}

\section{Computing initial guesses from Stokes expansion.}
\label{sec:stokes}

\begin{equation}
u_t+\left[f(u)\right]_{x} + \mathcal{L}u_x=0, \label{gen}
\end{equation}
	 The goal of this section is to explain how the idea of Stokes's approximation works in providing the initial data (guess) on wave and phase velocity for solving the equation~\eqref{gen} numerically.

	We will consider $\Lop$ being linear and self-adjoint Fourier multiplier operator, and a function $f$ that has degree of zeros $p \geq 2$:
\begin{align}
\widehat{\Lop u}(k) &= \alpha(k) \cdot \widehat{u}(k),\\
\langle\Lop u, v \rangle_{L^2(0,L)} &= \langle u, \Lop v\rangle_{L^2(0,L)}, \\
p &= \min\setc{k \in \NN}{f^{(k-1)}(0) =0 \text{ and } f^{(k)}(0)\neq 0}
\end{align}
	 Consider the equation~\eqref{gen} and its solution in the form $u(x,t) = \phi(x-ct)$, which is a traveling wave solution. Inserting $\phi(x-ct)$ into~\eqref{gen} leads to the equation:
\begin{equation}
-c\phi' + f'(\phi)\phi'+\mathcal{L}\phi'= 0,
\end{equation}
	 which can be integrated to give:
\begin{equation}
-c\phi + f(\phi)+\mathcal{L}\phi= B, \quad B = \mathrm{const}. \label{int_gen-B}
\end{equation}
	Consider $B = 0$ in equation~\eqref{int_gen-B}, and expansions of $\phi$ and $c$:
\begin{align}
\phi = \xi = \varepsilon\xi_1 + \varepsilon^2\xi_2+ \cdots, \label{fi_exp}\\
c = c_0 + \varepsilon c_1 + \varepsilon^2 c_2+\cdots. \label{c_exp}
\end{align}
	 The next step is to insert~\eqref{fi_exp} and~\eqref{c_exp} to the equation~\eqref{int_gen-B} and write out the terms at powers of $\varepsilon$. 
	The function $f(\phi)$ is expanded around zero and, therefore, will appear only in $\varepsilon^p$ terms. 
	Thus, the term at the first power of $\varepsilon$ reads: 
\begin{align}
\varepsilon: \quad -c_0\xi_1 + \mathcal{L}\xi_1 = 0, \label{epsilon1}
\end{align}
	 Hence, $c_0$ is an eigenvalue of the operator $\Lop$, regarded as defined on $L$-periodic functions. 
	 Taking the Fourier transform of the equation~\eqref{epsilon1} gives:
\begin{equation}
-c_0\widehat{\xi_1}(k) + \alpha(k)\widehat{\xi_1}(k) = 0 . \label{eps1-transformed}
\end{equation}
	 The equation~\eqref{eps1-transformed} has two trivial solutions: either $\xi_1(k) \equiv 0$ or $\alpha(k) \equiv c_0$. If we assume non-trivial $\xi_1$ and $\alpha(k) \neq \text{const}$, the following solves the problem:
\begin{align}
  \widehat{\xi_1}(k) = 2\pi\delta(k - k_0), \quad \text{and} \quad c_0 = \alpha(k_0), \label{solution-xi}
\end{align}
	 for some $k_0 \in \R$. 
	Since $\xi_1$ is the first-order approximation to $\phi$, the corresponding wave number should be equal to 1. The $L$-periodicity condition entails that $k_0 = 2\pi/L\cdot 1$. The spacial variable $x$ has to be scaled to $x'= L/2\pi x$, accordingly. From the solutions in~\eqref{solution-xi} we have 
\begin{equation}
	\xi_1(x') = e^{ik_0x'} = \cos(k_0x') + \sin(k_0x').
\end{equation}	
Considering the projection onto the space $\SN$ 
%
%
, we are led to choose $\xi_1(x') = \cos(k_0x')$.
%
%
%


	 For further analysis, let us define an operator $\Aop$ 
\begin{equation*}
	\Aop :=  - c_0 \mathcal{E} + \Lop,
\end{equation*}
where $\mathcal{E}$ is the identity operator. The operator $\Aop$ inherits the property of being 
self-adjoint from $\Lop$. Moreover, it follows from~\eqref{epsilon1} that $\Aop \xi_1 = 0$ and $\xi_1 \in \ker(\Aop)$. 
If $p>2$ then $f''(0)=0$ and the terms at $\epsilon^2$ are:
\begin{equation}
\Aop \xi_2 - c_1\xi_1 = 0.
\end{equation}

	 Taking scalar multiplication of the latter with $\xi_1$, one obtains:
\begin{align}
	&\langle\xi_1, \Aop\xi_2 \rangle_{L^2(0,L)} = c_1 \| \xi_1\|^2_{L^2(0,L)}, \\
	&\langle\xi_1, \Aop\xi_2\rangle_{L^2(0,L)} = \langle\Aop \xi_1, \xi_2\rangle_{L^2(0,L)} = \langle0, \xi_2\rangle_{L^2(0,L)} = 0.
\end{align}
	 As a result, one has $c_1 \| \xi_1\|^2_{L^2(0,L)} = 0$ and, hence, $c_1 = 0$.
	Repeating the same argument, it becomes clear that $c_k = 0$ for any $k \leq p-1$.
	Besides that, $\xi_2$ is in the kernel of $\Aop$, so it may be assumed to be proportional to $\xi_1$.

	The terms at order $\varepsilon^p$ are:
\begin{align}
  \Aop \xi_p - c_{p-1}\xi_1 + \frac{f^{(p)}(0)}{p!} \xi_1^p = 0. \label{epsp}
\end{align}
	 Let us denote for brevity
\begin{align*}
	f_p := \frac{f^{(p)}(0)}{p!}.
\end{align*}
	 Pairing~\eqref{epsp} with $\xi_1$ (and assuming $\| \xi_1 \|_{L^2(0,L)} = 1$) gives
\begin{align}
	c_{p-1} = f_p	\cdot \langle\xi_1^p, \xi_1\rangle_{L^2(0,L)}, \label{cp-1}
\end{align}
	 which gives us the value of $c_{p-1}$.
	It only remains to solve the following problem numerically in order to obtain $\xi_p$:
\begin{align}
  \Aop \xi_p = c_{p-1}\xi_1 - f_p\xi_1^p \label{epsilonp}
\end{align}

	For the last equation to be solved, the operator $\Aop$ has to be invertible. It is also required that $\langle\xi_1, \xi_p\rangle_{L^2(0,L)}=0$. 
	Therefore the solution is sought in the space orthogonal to $\ker(\Aop)$.
	Since $\Aop$ is still a Fourier multiplier operator, one can take the Fourier transform of the equation~\eqref{epsilonp} to find
\begin{align}
	\widehat\Aop(k) \widehat\xi_p(k)& = c_{p-1}\widehat\xi_1(k) - f_p				
	\widehat{\xi_1^p}(k), \\
	\widehat\xi_p(k)& = \widehat\Aop(k)^{-1} \left(c_{p-1}\widehat\xi_1(k) - f_p	
	\widehat{\xi_1^p}(k)\right).
\end{align}
	 Taking the inverse Fourier transform of $\widehat\xi_p(k)$ gives $\xi_p$.
	Since only even solutions of the problem are considered the cosine part of the Fourier transforms will be required.

	It is sufficient to use $\xi_1$ and $c_0$ as the initial guesses for the Newton method. However, it should be noted that for different values of $p$ the pair of parameters $\xi_p$ and $c_{p-1}$ are computed in different ways.
\begin{enumerate}[label=\alph*)]
	\item If $p = 2$, then $\xi_2$ is computed from~\eqref{epsilonp}, but $c_{p-1}$ here becomes zero. Therefore one has to consider the next level of the expansion $\varepsilon^p$.
	\item For odd values of $p$ the parameter $c_{p-1}$ can be computed from~\eqref{cp-1} and $\xi_p$ from~\eqref{epsilonp}.
	\item For even values of $p \geq 4$ the parameter $\xi_p$ can be computed, but $c_{p-1}$ may not be non-zero in general. 
	In such cases a different strategy of fixing the initial guess should be used.
\end{enumerate}


\section{Presentation of \textsf{SpecTraVVave} and its workflow}
\label{sec:implementation}
 
\subsection{Overview}
\begin{figure}[ht]
	\centering
\includegraphics[width=0.8\textwidth]{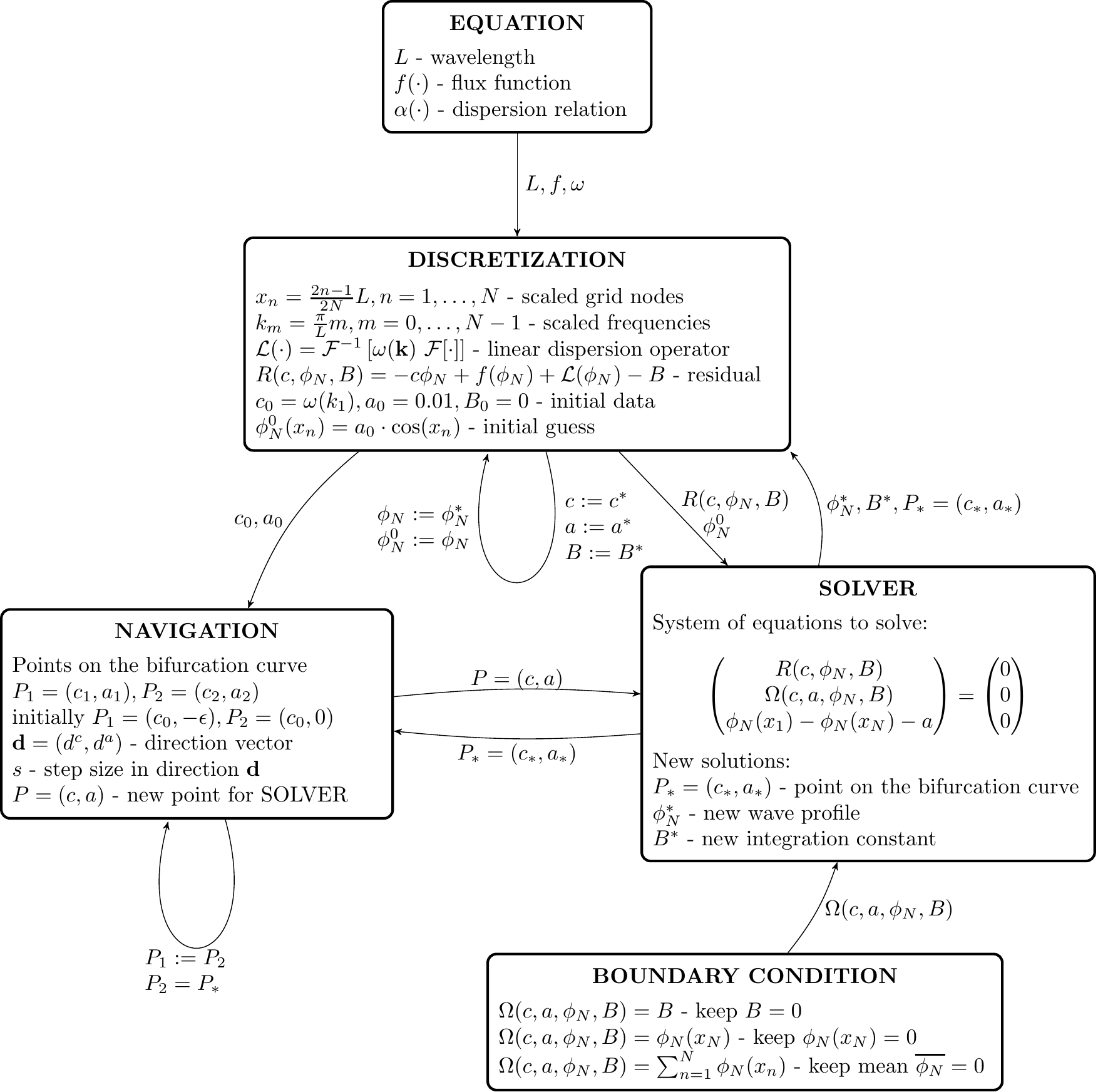} 
\caption{\small Overview of the \textsf{SpecTraVVave} package.}     
\label{schematic}
\end{figure}

There are several classes in the \textsf{SpectraVVave} package. An overview of the program is shown 
in \autoref{schematic}. The workflow begins with defining a flux function $f$ and the Fourier multiplier function $\alpha$ to set up an equation. 
The traveling wave solution is characterized by the wavelength $L$ 
and a boundary condition $\Omega(c, a, \phi_N, B)$. These parameters are fixed for a given problem.
The defined equation is then discretized. 
The \pyth{Discretization} object contains all required elements such as grid points, wave-numbers and the  discrete linear operator.
		
The initial guess and the equation's residual are passed from the \pyth{Discretization} to the \pyth{Solver} object. 
The \pyth{Navigation} object is responsible for finding good initial guesses for $c$ and $a$ that are  
passed to the \pyth{Solver} object. 
The \pyth{Solver} object applies Newton's method to find a solution to the system of equations~\eqref{system-N3}. 
	
The new solution is sent back to the  \pyth{Discretization} and \pyth{Navigation} objects, where variables get updated.
All computed solutions are stored for further analysis. 
This finishes one iteration. 
For the next iteration the updated variables are used and a new solution is found. 
The process may be continued as long as the Jacobian of the problem is non-singular. 
 
 \subsection{Class Description}

We present an overview of the classes used in \textsf{SpectraVVave} package. Note that, since the package is under continuous modification and development, we describe here only the  basic classes and functions the package. 
We refer to the package repository \cite{github} for up-do-date tutorials and installation instructions.
 
 The \pyth{Equation} class is the general class for all model equations. 
Its only role is to store a parameter $L$, the wavelength:

\begin{python}
class Equation(object):
    def __init__(self, L):
        self.length = L
\end{python}
A subclass of the \pyth{Equation} class has to implement two functions, \pyth{compute_kernel} and \pyth{flux}. 

The KdV model equation
\begin{align}
	u_t + \frac{3}{2} u u_x + u_x + \frac{1}{6} u_{xxx} = 0
\end{align}
with $f(u) = \frac{3}{4} u^2$ and $\widehat{\Lop u}(k) = (1 - \frac{1}{6} k^2)\widehat{u}(k)$ is presented in the program as a subclass of the \pyth{Equation} class:
\begin{python}
class KDV(Equation):
    def compute_kernel(self, k):
        return 1-1/6*k*k

    def flux(self, u):
        return 3/4*u*u  
\end{python}

On can then create an object of the class \pyth{KDV} with the command:
\begin{python}
kdv = KDV(L=np.pi)
\end{python} 

The solver will compute only a half of a solutions profile. The fundamental wavelength of the solutions of the defined equation will be equal to $2\pi$. 

In order to find solutions with specific features, boundary conditions are introduced as separate classes.
For instance, the boundary condition specifying a constant of integration is implemented as follows:
\begin{python}
class Const(object):
""" The boundary condition under which the constant of integration 
(B) is always set to zero. """
    def __init__(self, level=0):
        self.level = level

    def enforce(self, wave, variables, parameters):
    """ Enforces the Const boundary condition. """
        return np.hstack([variables[0] - self.level])

    def variables_num(self):
    """ The number of additional variables that are required to construct 
    the Const boundary condition. """
        return 1
\end{python}
A \pyth{Const} boundary condition object is created as follows:
\begin{python}
boundary_condition = Const()
\end{python}

The next step is to create an object of \pyth{Discretization} class, which is initialized with a model equation such as \pyth{kdv_model} and the number of grid points. The main parts of the class are the following:
\begin{python}
class Discretization(object):
    def __init__(self, model_equation, grid_size):
        self.equation = model_equation
        self.size = grid_size

    def operator(self, u):
        u_ = scipy.fftpack.dct(u, norm='ortho')
        Lv = self.fourier_multiplier()*u_
        result = scipy.fftpack.idct(Lv, norm='ortho')
        return result

    def residual(self, u, wavespeed, const_B):
        residual = - wavespeed*u + self.equation.flux(u) 
                                        + self.operator(u) - const_B
        return  residual
\end{python}
The call \pyth{Discretization.operator(u)} computes $\Lop u $ as the inverse transform of a transformed convolution 
\begin{equation}
\Lop u = \mathcal{F}^{-1}[ \widehat{\Lop u}(k)] = \mathcal{F}^{-1}\left[ \alpha(k)\cdot\widehat{u}(k)\right].
\end{equation}
The result of the call \pyth{Discretization.residual(u, wavespeed, const_B)} is then used in the \pyth{Solver} class. 
An object of the \pyth{Solver} class is initialized with an object of the \pyth{Discretization} class, and a boundary condition object.

\begin{python}
class Solver(object):
    def __init__(self, discrete_problem, boundary_condition):
        self.discretization = discrete_problem
        self.boundary = boundary_condition

    def solve(self, guess_wave, parameter_anchor, direction):
    """ Runs a Newton solver on a system of nonlinear equations once. 
    Takes the residual(vector) as the system to solve. parameter_anchor 
    is the initial guess for (c,a) values and it is taken from the 
    Navigation class. """
        size = len(guess_wave)
        self.discretization.size = size

        def residual(vector):
        """ Contructs a system of nonlinear equations. First part, 
        main_residual, is from given wave equation; second part, 
        boundary_residual, comes from the chosen boundary conditions. """
        . . .   
        return np.hstack([main_residual, boundary_residual, 
                                          amplitude_residual])  
    . . .
    return new_wave, new_boundary_variables, new_parameter
\end{python}

Some omitted parts in the above script are substituted by \textsf{'. . .'} sign. Each iteration on a \pyth{Solver} object is run from a \pyth{Navigation} object, which takes the \pyth{Solver} object for initialization. 
%
\begin{python}
class Navigation(object):
""" Runs the Solver and stores the results. """

    def __init__(self, solver_object, size=32):
    self.solve = solver_object.solve # function to run Newton method
    self.size = size # size for navigation
    . . . 
    def run_solver(self, current_wave, pstar, direction):
        new_wave, variables, p3 = self.solve(current_wave, pstar, direction)
        return new_wave, variables, p3
\end{python}

All the above classes can be modified and developed further, new classes may be defined as well.

\subsection{Detailed Workflow}

The workflow with the package consists of three basic steps:
\begin{enumerate}
\item Once the necessary classes have been imported in the current namespace, generate all necessary objects:
\begin{python}
equation = KDV(L=np.pi)
boundary_condition = Const()
discretization = Discretization(equation, grid_size=64)
solver = Solver(discretization, boundary_condition)
navigator = Navigation(solver)    
\end{python}
\item Choose a number of iterations, i.e., the number of solutions to compute, and run the solver:
\begin{python}
n_iter = 50
navigator.run(n_iter)
\end{python}
\item All computed solutions are stored in \pyth{navigation_object}
\begin{python}
last_computed = -1    
wave_profile = navigator[last_computed]['solution']
wave_speed = navigator[last_computed]['current'][0]
wave_amplitude = navigator[last_computed]['current'][1]
\end{python}
\end{enumerate} 

We further refer to the code repository \url{https://github.com/olivierverdier/SpecTraVVave} for up-to-date instructions on how to run the code.




\end{document}